\newfont{\Bbb}{msbm10 scaled\magstep1}

\newcommand{\CC}{{\mbox{\Bbb C}}}

\newcommand{\PP}{{\mbox{\Bbb P}}}

\newcommand{\NN}{{\mbox{\Bbb N}}}
\newcommand{\QQ}{{\mbox{\Bbb Q}}}

\newcommand{\ZZ}{{\mbox{\Bbb Z}}}

\def\kk{{\bf k}}
\def\TT{{\bf T}}

\def\vv{{\bf v}}
\def\bt{{\bf t}}
\def\dM{{M^\vee}}

\def\vv{{\bf v}}
\def\ulamb{{\underline{\lambda}}}
\def\umu{{\underline{\mu}}}
\def\unu{{\underline{\nu}}}
\def\lra{\longrightarrow}

\def\ee{\epsilon}
\def\proof{{\bf Proof.}\,\,}
\def\qed{\hfill{{\bf QED}}}

\def\be{\begin{equation}}
\def\ee{\end{equation}}
\def\kformep{\epsilon^1\wedge\ldots\wedge\epsilon^k}

\def\ikformep{\epsilon^{i_1}\wedge\ldots\wedge\epsilon^{i_k}}
\def\jkformep{\epsilon^{j_1}\wedge\ldots\wedge\epsilon^{j_k}}

\def\ihikformep{\epsilon^{i_1+h_1}\wedge\ldots\wedge\epsilon^{i_k+h_k}}
\def\rikformep{\epsilon^{1+r_1}\wedge\ldots\wedge\epsilon^{k+r_k}}

\def\wM{\bigwedge M}

\def\wkM{\bigwedge^k M}
\def\wkdM{\bigwedge^k M^\vee}

\def\ovW{\overline{W}}

\documentclass[12]{article}
\usepackage{graphicx}
\usepackage{eufrak}

\newfont{\eulerfraktur}{eufm10 scaled\magstep1}

\newtheorem{thm}{Theorem}[section]
\newtheorem{defin}{Definition}[section]
\newtheorem{example}{Example}[subsection]
\newtheorem{prop}{Proposition}[section]
\newtheorem{corol}{Corollary}[section]
\newtheorem{lemma}{Lemma}[section]
\newtheorem{rmk}{Remark}[section]

\begin{document}
\title{\normalsize{\Large The Algebra of Schubert Calculus}
}
\author{{\small LETTERIO GATTO}\thanks{Work partially
sponsored by MIUR  (Progetto Nazionale ``Geometria sulle Variet\`a
Algebriche",  coordinatore Sandro Verra),
and supported by GNSAGA-INDAM.}\\{ }
\\ Dipartimento di Matematica,
Politecnico di Torino,
\\
{\normalsize Corso Duca degli Abruzzi 24, 10129 Torino -- (ITALY)}\\
}
\date{}

\maketitle
\begin{abstract}
\noindent
A flexible unified framework for
both classical and quantum Schubert calculus is proposed. It is based on a
natural combinatorial approach relying on the {\em Hasse-Schmidt} 
extension of a
certain family of pairwise commuting endomorphisms of an infinite free
$\ZZ$-module $M$ to its exterior algebra $\wM$.
\end{abstract}

\section{Introduction}

\subsection{The Goal} Let $V$ be an $n$-dimensional vectorspace over an
algebraically closed field $\kk$ and let
$G:=G_k(V)$ be the grassmannian variety
parametrizing
$k$-dimensional subspaces of it.
The
main goal of this paper is to provide a flexible unified framework for
both classical and quantum Schubert calculus via a new combinatorial approach
based on elementary considerations of linear algebra.

To give the flavour of the results, our favourite example is that of an
infinite
$\ZZ$-module
$M$, freely spanned by $(\epsilon^1,\epsilon^2,\epsilon^3,\ldots)$, equipped
with the {\em shift} endomorphism $D_1:M\lra M$, defined on generators by
$D_1(\epsilon^i)=\epsilon^{i+1}$.  If one extends it to $\bigwedge^2M$, by
imposing Leibniz's rule with respect to $\wedge$, one has, for instance:
\be
D_1(\epsilon^2\wedge\epsilon^4)=D_1\epsilon^2\wedge\epsilon^4+\epsilon^2\wedge
D_1\epsilon^4=
1\cdot\epsilon^3\wedge\epsilon^4+1\cdot\epsilon^2\wedge\epsilon^5.
\label{eq:forprima}
\ee
In spite of computational simplicity,
the coefficients occurring on the right hand side of
equality~(\ref{eq:forprima}), have a remarkable geometrical meaning. In
fact,
let \linebreak $(P_1,P_2,\ell_1,\ell_2,\ell_3)$ be a sufficiently general
configuration of two points and three lines in $\PP^3$  (e.g. the $\ell_i$'s
and the line through $P_1$ and $P_2$ are pairwise skews): then the
$1$ multiplying $\epsilon^3\wedge\epsilon^4$ is the degree of the intersection
of the Schubert variety
$\ovW_1(P_1)$, closure in $G_1(\PP^3_{\bf C})$ of all the lines of
$\PP^3_{\bf C}$ passing through $P_1$, and of the Schubert variety
$\ovW_{(2,1)}(P_2,\ell_1,\ell_2)$, closure of all the lines of
$\PP^3$ meeting $\ell_1$ and $\ell_2$ and passing through $P_3$; furthermore,
the
$1$ multiplying $\epsilon^2\wedge\epsilon^5$ is the number of degree $1$
rational maps
$f:\PP^1_{\bf C}\lra G_1(\PP^3_{\bf C})$, such that
$f(0)\in \ovW_1(P_1)$, $f(1)\in \ovW_{(2,1)}(P_2,\ell_1,\ell_2)$ and
$f(\infty)= [\ell_3]\in G_1(\PP^3_{\bf C})$.

The reason why is so, as will be explained in the paper, is that the $k$-th
exterior power of a free $\ZZ$-module $M_n$ of rank
$n$ is a free principal module over the Chow intersection ring $A^*(G)$ of
$G$, the latter operating on the former as a ring of differential
operators. A similar statement holds for $\wkM\otimes_{\bf Z}\ZZ[q]$ ($q$ an
indeterminate over $\ZZ$), which turns out to be a principal free module over
the small quantum intersection ring $QA^*(G)$ of $G$. Although  these
facts are rather elementary (or perhaps for this reason), they seem to have
been at least ignored in the literature, if not unknown. By contrast, we
contend that they summarize the essential algebraic content of Schubert
calculus for grassmannians, in its classical and quantum formulation.

\subsection{The Algebra of \ldots} Following Lang~(\cite{Lang}, p.~556 ) an
$A$-module
$M$ is said to be {\em principal}, generated by $m_0\in M$, if the natural map
  $E_{m_0}:A\lra M$, sending $a\in A$
onto $am_0$, is an epimorphism. Any principal $A$-module
inherits therefore a structure of commutative ring induced by the isomorphism
$\Pi:A/\ker E_{m_0}\lra M$, by setting $m_1*m_2=\Pi(\Pi^{-1}(m_1)\cdot
\Pi^{-1}(m_2))$. The inspiring idea consists in interpreting the 
isomorphism $\Pi$ as
if it were  Poincar\'e duality between homology and cohomology of some smooth
algebraic variety.

Within this philosophy, the main actor of our story is a $\ZZ$-module
$
M=\bigoplus_{i\geq 1}\ZZ\epsilon^i$ together with a $\ZZ$-module
homomorphism
\[
D_t=\sum_{i\geq 0}D_it^i: M\lra M[[t]],
\]
where each $D_i\in End_{\bf Z}(M)$
is defined on generators by $D_i\epsilon^j=\epsilon^{i+j}$. Clearly
$D_i=D_1^i$ and it is easy to show that there is a  $\ZZ$-module
isomorphism\linebreak
$E_{\epsilon^1}:\ZZ[D_1]\lra M$, given by $P(D_1)\mapsto P(D_1)\epsilon^1$.
One says that the data $(M,D_t)$ is an {\em intersection pair} and that
${\cal A}^*(M,D_t):=\ZZ[D_1]$, is its {\em intersection ring}.

Next, and this is the key point, one proceed to extend such a
$D_t$ to the exterior algebra $\wM$ of $M$, by making it into a
{\em Hasse-Schmidt derivation} (Cf.~\cite{Mats}, p.~208, or the 
original paper~\cite{HasSm}, for commutative algebras). In other 
words, one requires that
$D_t(\alpha\wedge\beta)=D_t(\alpha)\wedge D_t(\beta)$, for all $\alpha,\beta\in
\wM$. Such an extension of $D_t$ will be called {\em Schubert derivation} on
$\wM$:  {\em derivation} because the components
$D_i$ of
$D_t$ behave as $i$-th derivatives do:
\be
D_i(\alpha\wedge\beta)=\sum_{\matrix{i_1+i_2=i\cr i_1,i_2\geq
0}}D_{i_1}\alpha\wedge D_{i_2}\beta ,\label{eq:intrler}
\ee
and {\em Schubert}, because Leibniz's
rule~(\ref{eq:intrler}) implies {\em Pieri's formula} of both classical and
quantum Schubert calculus -- see below.
Therefore $\wkM$ turns out to be  a module over
$\ZZ[D]:=\ZZ[D_1,D_2,\ldots]$, the commutative $\ZZ$-subalgebra of $End_{\bf
Z}(\wkM)$ generated by the $D_i$'s. The key point is that $(\wkM,D_t)$ is
an {\em intersection pair}, too, i.e. $\wkM$ is a principal
$\ZZ[D]$-module generated by\linebreak
$\kformep$. This is  a consequence of

\medskip
\noindent
{\bf Theorem~\ref{giamth}} {\em {\em Giambelli's formula} holds:
\[
\rikformep=\Delta_{(r_k\ldots r_1)}(D)\cdot \kformep,
\]
where \[
\Delta_{(r_k\ldots r_1)}(D)=\left|\matrix{D_{r_1}&D_{r_2+1}&\ldots&
D_{r_k+k-1}\cr D_{r_1-1}&D_{r_2}&\ldots& D_{r_k+k-2}\cr
\vdots&\vdots&\ddots&\vdots\cr
D_{r_1-k+1}&D_{r_2-k+2}&\ldots&D_{r_k}}\right|.
\]
is  {\em Giambelli's determinant}.
}

\medskip
\noindent
One then proves that the intersection ring ${\cal A}^*(\wkM,D_t) =\ZZ[D]/\ker
(E_{\kformep})$ of the pair $(\wkM,D_t)$ is  isomorphic
(Proposition~\ref{propgir}) to the polynomial ring $\ZZ[D_1,\ldots,D_k]$,
since, for each $h>k$,
$D_h$ admits an explicit polynomial expression $D_h=D_h(D_1,\ldots, D_k)$ of
$D_1,\ldots,D_k$ only.

Once one is given of such a general framework, it is
natural to look at finite rank situations: two more intersection pairs can
  in fact be costructed out of
$(\wkM,D_t)$. They will be denoted
$(\wkM_n,D_t)$ and
$(\wkM[D_n],D_t)$ respectively, where $M_n:=M/D_nM$  and $M[D_n]$ is $M$
itself thought of as a $\ZZ[D_n]$-module.

That $\wkM_n$ is  a principal
$\ZZ[D_1,\ldots,D_k]-$module generated by
$\kformep$ is a consequence of the isomorphism
$E_{\kformep}:\ZZ[D_1,\ldots,D_k]\lra \wkM$ composed with the canonical
projection
$\wkM\lra
\wkM_n$. More precisely:

{\bf Proposition~\ref{presip}} {\em The intersection ring of 
$(\wkM_n,D_t)$ is:
\be
{\cal A}^*(\wkM_n,D_t)={\ZZ[D_1,\ldots,D_k]\over(D_{n-k+1},\ldots,
D_n)}.\label{eq:inpresip}
\ee
}

A few words about the proof. To check that, for each $i\geq 0$, the polynomial
$D_{n-k+i}$ belongs to the ideal of  relations is an easy matter: on one hand
$\wkM_n$ is isomorphic to the quotient of
$\wM$ by the ideal
$\bigwedge^{k-1}M\wedge D_nM$ (generated by all $k$-vectors $\ikformep$ such
that $i_k>n$) and, on the other hand:
\begin{eqnarray*}
D_{n-k+i}\kformep&=&\epsilon^1\wedge\ldots\wedge\epsilon^{k-1}\wedge
\epsilon^{n+i}=\\ &=&\epsilon^1\wedge\ldots\wedge\epsilon^{k-1}\wedge
D_n\epsilon^{i}\in \bigwedge^{k-1}M\wedge D_nM.
\end{eqnarray*}
To show that, indeed, $(D_{n-k+1},\ldots,D_n)$ is a complete set of relations,
one exploits instead the very shape of Giambelli's determinant
(Prop.~\ref{propideale}).

As for $M[D_n]$, it is a free $\ZZ[D_n]$-module of rank
$n$ generated by $(\epsilon^1,\ldots,\epsilon^n)$ and
$(\wkM[D_n],D_t)$ is a $\ZZ[D_n]$-intersection pair, too. In fact $\wkM[D_n]$
is a principal
$\ZZ[D_n][D_1,\ldots,D_k]$-module generated by $\kformep$ and its intersection
ring is isomorphic to:
\be
\ZZ[D_n][D_1,\ldots,D_k]\over{(D_{n-k+1},\ldots, D_{n-1})}\label{eq:presdn}
\ee
The latter can also be written via renaming $D_n$ by an auxiliary
indeterminate $q$ over $\ZZ$, getting:
\be
{\cal A}^*(\wkM_n[q],D_t):={\ZZ[D_n][D_1,\ldots,D_k]\over{(D_{n-k+1},\ldots,
D_{n-1})}}={\ZZ[q][D_1,\ldots,D_k]\over{(D_{n-k+1},\ldots, D_{n-1},
D_n-q)}},\label{eq:preseq}
\ee
where  $M_n[q]\cong M[D_n]/(D_n-q)\cong M_n\otimes_{\bf Z}\ZZ[q]$.
Up to a few changes, similar arguments used in the ``classical" case, show that
$(D_{n-k+1},\ldots, D_{n-1})$ is indeed a complete set of relations for the
presentation~(\ref{eq:presdn}) of ${\cal A}^*(\wkM[D_n],D_t)$.

\subsection{\ldots Schubert Calculus}

Expressions~(\ref{eq:inpresip}) and~(\ref{eq:preseq}), up to renaming $D_i$ and
changing $q$ by $(-1)^{k-1}q$, are exactly the expressions
  of $A^*(G)$ and $QA^*(G)$ respectively (Corollaries~\ref{corisoclas}
and~\ref{corwstder}). However, the relationship of the above algebraic model
with Schubert calculus is not  purely formal. As a matter of fact the
$k$-vectors $\{\ikformep\}_{1\leq i_1<\ldots<i_k\leq n}$ corresponds to
Schubert cycles. More precisely, recall that the Chow group
$A_*(G)$ is freely generated over the integers by the classes modulo rational
equivalence of  Schubert varieties (\cite{Fu1}, p.~27 and the references
therein). These are parametrized by partitions of lenght
$\leq k$ and weight
$\leq k(n-k)$, i.e. by non increasing sequences of $k$ nonnegative integers
$\ulamb=(\lambda_1\geq
\ldots\geq\lambda_k)$, such that $0\leq |\ulamb|=\lambda_1+\ldots+\lambda_k\leq
k(n-k)$. The Schubert variety $\overline{W_\ulamb({\cal E})}$  is the closure
of the locus of
$k$-planes verifying certain incidence conditions with respect to  a given flag
\be
{\cal E}:={\rm span}\{e_1,e_2,\ldots,e_n\}\supset
{\rm span}\{e_2,\ldots,e_n\}\supset\ldots\supset {\rm span}\{e_n\}\supset
{\rm span}\{0\},\label{eq:intradapt}
\ee
prescribed by the partition $\ulamb$ itself. Let $\sigma_\ulamb\in A^*(G)$ be
the Schubert (co-)cycle such that $\sigma_\ulamb\cap 
[G]=\ovW_\ulamb$. The intersection ring $A^*(G)$ (resp.
the quantum intersection ring
$QA^*(G)$) of the grassmannian is generated over the integers (resp. over
$\ZZ[q]$) by the {\em special Schubert cycles} $\sigma_i:=\sigma_{(i0\ldots
0)}$
and  is isomorphic to the ring~(\ref{eq:inpresip}) (resp. to the
ring~(\ref{eq:preseq})) via $\sigma_i\mapsto D_i$ (resp. $\sigma_i\mapsto
D_i$ and $q\mapsto (-1)^{k-1}q$), so proving what announced: the 
groups $\wkM_n$
and
$\wkM_n[q]$ are principal modules over $A^*(G)$ and $QA^*(G)$
respectively, with these operating as
rings of derivations. However, the bridge with Schubert calculus is due
to the very geometrical interpretation of the $k$-vectors $\ikformep$, shortly
described below.

Let
$(e_1,\ldots, e_n)$ be an adapted basis of  a complete flag
${\cal E}$ of subspaces of $V$   as in~(\ref{eq:intradapt}) and let
$(\epsilon^1,\ldots,\epsilon^n)$ be the dual basis. Let $M_n$ be the integral
lattice generated by the
$(\epsilon^j)$. There is then a canonical isomorphism
\[
[\,\,]:\wkM_n\lra A_*(G)
\]
associating to each $\ikformep$ a set of equations
defining scheme theoretically the Schubert variety $\ovW_\ulamb$, with
$\ulamb=(i_k-k,\ldots, i_1-1)$. Therefore $A_*(G)$ (resp.
$QA^*(G)$) comes equipped with two  structures of
$A^*(G)$-module (resp. $QA^*(G)$-module): the former induced by the ({\em
cap}-)intersection product $\cap$, the latter gotten by setting, for each
$P(D)\in{\cal A}^*(\wkM_n,D_t)\cong A^*(G)$:
\[
P(D)\cdot\ovW_\ulamb=[P(D)\ikformep].
\]

One is then finally left with showing that the two module structures
coincide, proving that the algebra of Schubert derivations on $\wkM_n$ (resp.
$\wkM_n[q]$) is indeed Schubert calculus (resp. quantum Schubert calculus).
Clearly the proof amounts to show that:
\be
\sigma_h\cap[\ikformep]=D_h[\ikformep]=[D_h(\ikformep)],\label{eq:equivmod}
\ee
because any $\sigma_\ulamb$ is a polynomial expressions in the special Schubert
cycles.

Left hand side of formula~(\ref{eq:equivmod}) can be computed via Pieri's
formula for Schubert calculus. Therefore, proving equality~(\ref{eq:equivmod})
is equivalent to prove a {\em Pieri's formula for Schubert 
derivations}. It is a
consequence of an elementary combinatorial observation: when computing
$D_h(\ikformep)$, cancelations of summands may occur, due to the
$\ZZ_2$-symmetry of $\wM$. Pieri's formula for ${\cal S}$-derivations, instead,
compute $D_h(\ikformep)$ by listing the surviving summands only. This makes
possible the formal proof that, indeed:
\[
\sigma_h\cap \ovW_\ulamb=D_h[\ikformep],
\]
where $\ulamb=(i_k-k,\ldots, i_1-1)$, concluding the argument. The same
can be argued with very slight modifications in the case of quantum Schubert
calculus, where Bertram-Pieri's formula simply translates in computing a
derivative of a $k$-vector $\ikformep\in\wkM$ and  replacing
$\epsilon^{nj+i}$ ($1\leq i<n$) with
$(-1)^{k-1}q^j\epsilon^i$ in any summand
of the output.

\subsection{ } The paper is organized as follows. Section~\ref{Sect2} is
devoted to preliminaries and notation, including a quick review of (quantum)
Schubert calculus. Sect.~\ref{Sect3} is the core of the paper:
the {\em canonical intersection pair} is studied, together with its
Hasse-Schmidt extension. Sect.~\ref{Sect4} studies the intersection pair $(M_n,
D_t)$ and its quantum deformation $(M_n[q],D_t)$ while Sect.~\ref{Sect5}
formally establishes the already evident bridge between Schubert 
derivations and
Schubert calculus. The paper is concluded with a couple of examples showing the
application of the methods.

\medskip
\noindent
{\large\bf Acknowledgment.} The author is 
grateful to Prof.~Dan Laksov for his friendly support and for many 
precious advises and criticism, which substantially improved the 
shape of the paper.

\section{Preliminaries and Notation}\label{Sect2}

\subsection{Partitions}\label{sect1.1}

In this paper we shall deal with partitions, denoted by underlined small
greek letters ($\ulamb$, $\umu$, $\unu$,\ldots), and with finite increasing
sequences of positive integers, denoted by capital roman letters ($I$, $J$,
$K$,\ldots), and shortly named {\em symbols}.
  A partition
$\ulamb$ is a sequence of non negative integers
$(\lambda_1\geq\ldots\geq\lambda_k\geq\ldots)$ such that all but
finitely many terms are equal to $0$. The {\em lenght} of a partition is
$\ell(\ulamb)=\sharp\{i:\lambda_i\neq 0\}$, its {\em weight}  is
$|\ulamb|=\sum\lambda_i$.  Let ${\cal P}$ be the
set of all partitions. The following notation shall be used:
\[
{\cal P}_h=\{\ulamb\in {\cal P}:|\ulamb|=h\}\,\,(h\geq 0);\quad {\cal
L}_k=\{\ulamb\in{\cal P}: \ell(\ulamb)\leq k\}.
\]
If $\ulamb\in{\cal L}_{\leq k}$, set
$r_i(\ulamb)=\lambda_{k+1-i}$ ($1\leq i\leq k$), in order to write it
either as
$\ulamb=(\lambda_1\lambda_2\ldots\lambda_k)$ or as $(r_kr_{k-1}\ldots r_1)$
(skipping  ``$\ulamb$" from the notation),by possibly adding a final string
of zeros if $\ell(\ulamb)<k$.   The
$\lambda_i$'s or the $r_i$'s are said to be the {\em parts}.
Let $\ulamb=(r_k\ldots r_1)\in{\cal L}_{\leq k}$. The symbol $I(\ulamb)=(1\leq
i_1\leq \ldots\leq i_k)$ associated to $\ulamb$ is defined by $i_j=r_j+j$.
Conversely, given any symbol
$I=(1\leq i_1< \ldots< i_k)$ one associates a partition
$\ulamb(I)=(r_k,\ldots,r_1)$, where $r_j=i_j-j$. The {\em weight} of a symbol
is
\[
wt(I)=|\ulamb(I)|=r_1+\ldots+r_k=(i_1-1)+\ldots+(i_k-k)=|I|-{k(k-1)\over 2},
\]
where $|I|=i_1+\ldots+i_k$.
When needed a
partition $\ulamb\in{\cal L}_{\leq k}$ shall also be written as
$(0^{m_0}1^{m_1}2^{m_2}...)$, where
$m_i=\sharp\{j:\lambda_j=i\}$ is the {\em multiplicity} which the integer $i$
occurs with in $\ulamb$.  For
instance, if $\ulamb=(444311100)$, one may
also write $(0^21^{3}3^14^3)$.
The partition having
$k$ parts equal to
$n$ can hence be written as $(n^k)$ instead of $(n,n,\ldots, n)$. A positive
integers shall be identified with a partition with just one part. The 
dual of a partition $\ulamb=(r_k,\ldots,r_1)$ is the
partition $\ulamb^\vee=(n-r_1,n-r_2\ldots,n-r_k)$.
As in
\cite{MacDon}, $\underline{\delta}_k\in {\cal L}_{\leq k}$ will stand
for the partition
$(k-1,k-2,\ldots,0)$.

Let $A$ be a commutative ring with unit and let $c:\ZZ\lra A$ be a
sequence.  Set $c_i=c(i)$.
To any partition of lenght $k$,  one may
associate the element
$\Delta_\ulamb(c)\in A$, defined as:
\be
\Delta_{\ulamb}(c)=\left|\matrix{c_{r_1}&c_{r_2+1}&
\ldots&c_{r_k+k-1}\cr
c_{r_1-1}&c_{r_2}&\ldots&c_{r_k+k-2}\cr
\vdots&\vdots&\ddots&\vdots\cr
c_{r_1-k}&c_{r_2-k+1}&\ldots&c_{r_k}}\right|.\label{eq:gdet}
\ee
A multi-index of lenght $k$ is any point of $\NN^k$. Multindices
can be partially ordered via the {\em Bruhat order}, by declaring that
$(i_1,\ldots,i_k)\preceq (j_1,\ldots,j_k)$ if and only if
$i_\alpha\leq j_\alpha$ for each $\alpha\in\{1,\ldots,k\}$. Similarly
$(i_1,\ldots,i_k)\prec (j_1,\ldots,j_k)$ if and only if
$(i_1,\ldots,i_k)\preceq (j_1,\ldots,j_k)$ and there is
$\alpha\in\{1,\ldots,k\}$ such that $i_\alpha<j_\alpha$.

\subsection{Symmetric Polynomials}  \label{sect2.2} Let us recall
some notation borrowed from \cite{MacDon}.  To each
$\ulamb\in {\cal L}_{\leq k}$ one  associates the following symmetric
polynomials in the set of indeterminates $X:=(x_1,\ldots, x_k)$.
\begin{enumerate}
\item The {\em monomial symmetric functions}
$m_\ulamb(X)=\sum x_1^{\eta(\lambda_1)}\cdot\ldots\cdot
x_k^{\eta(\lambda_k)}$, summed over all distinct
permutations
$\eta$ of $(\lambda_1,\ldots,\lambda_k)$ (\cite{MacDon}, p. 11);
\item The {\em complete symmetric polynomials}
$h_\ulamb(X)=h_{\lambda_1}(X)\cdot\ldots\cdot h_{\lambda_k}(X)$,
where\newline
$h_n(X):=\sum_{|\ulamb|=n}m_\ulamb(X)$ (\cite{MacDon}, p. 14).
\item The {\em Schur polynomials}:
\[
s_\ulamb={\Delta_{\lambda+\underline{\delta}}(X)\over
\Delta_{\underline{\delta}}(X)}
\]
where $\Delta_{\underline{\delta}}(X)=\prod_{i<j}(x_i-x_j)$ is the usual
Vandermonde determinant (\cite{MacDon}, p. 23 ff.).

\end{enumerate}
In the sequel we shall need the following results:
\begin{enumerate}
\item The polynomials $s_\ulamb$, $h_\ulamb$ and $m_\ulamb$ are a
$\ZZ$-basis of the symmetric polynomials in $X$.(\cite{MacDon}, p. 11, 15, 24)
\item The $s_\ulamb$ are related with the $h_i$ via the determinantal
formula: $s_\ulamb=\Delta_\ulamb(h)$. In particular $s_i:=s_{(i)}=h_i$
for each $i\geq 0$, (\cite{MacDon}, p. 25).
\end{enumerate}

\subsection{(Quantum) Schubert Calculus}
Let $G:=G_k(V)$ be the grassmannian variety parametrizing $k$-planes in a $n$
dimensional vectorspace $V$. It is a smooth homogeneous variety and its Chow
ring $A^*(G)$ is freely generated by the {\em Schubert (co)cycles}
$\{\sigma_\ulamb:\ulamb\in{\cal L}_{\leq k}\}$ defined by the equality:

\be
\sigma_\ulamb\cap [G]=\overline{W}_\ulamb\in A_*(G),\label{eq:shc1}
\ee
where $\ovW_\ulamb$ is the class modulo rational equivalence of a {\em Schubert
variety}, to be defined below.
Let
\be
{\cal E}: \quad V=E_0\supset E_1\supset\ldots\supset E_n=(0),\label{eq:flag}
\ee
be any complete flag of $V$ ( $i={\rm
codim}\,E_i$).

The flag~(\ref{eq:flag}) induces a chain of inequalities:
\[
k=\dim([\Lambda]\cap E_0)\geq \dim([\Lambda]\cap E_1)\geq\ldots\geq
\dim([\Lambda]\cap E_n)=0.\label{eq:flagdim}
\]
One says that $i_j$ is a (${\cal E}$-){\em Schubert jump} of 
$[\Lambda]\in G$ if
and only if
\[
\dim([\Lambda]\cap E_{i_j-1})>\dim([\Lambda]\cap E_{i_j}).
\]
Since the
dimension decreases stepwise no more than $1$, there are exactly $k$ jumps:
$
1\leq i_1<\ldots<i_k\leq n.
$
The sequence of the Schubert jumps of $[\Lambda]$,
$I_{\cal E}([\Lambda])=(i_1,\ldots,i_k)$, is
said to be the {\em Schubert symbol} of
$[\Lambda]\in G$.
The ${\cal E}$-Schubert variety
$\overline{W_\ulamb({\cal E})}$ is the closure in
$G$ of the {\em Schubert cell}
\be
W_\ulamb({\cal E})=\{[\Lambda]\in G:I_{\cal
E}([\Lambda])=I(\ulamb)\},
\ee
where $I(\ulamb)$, the {\em symbol} associated to
$\ulamb$ (Cf. Sect~\ref{sect1.1}), is equal to the Schubert symbol
$I_{\cal E}([\Lambda])$. In the sequel we shall also denote a Schubert cell
as $W_I$ meaning $\{[\Lambda]\in G:I_{\cal
E}([\Lambda])=I\}$,

The class modulo rational equivalence  of $\overline{W_\ulamb({\cal
E})}$ (the closure in $G$ of the ${\cal E}$-Schubert cell) does not depend on
the flag chosen, but only on
$\ulamb$: it  is said to be a {\em Schubert cycle} and shortly denoted by
$\ovW_\ulamb$. For any complete flag ${\cal E}$ of $V$ and any $[\Lambda]\in G$
there exists
$\ulamb\in {\cal L}_{\leq k}$ such that $[\Lambda]\in W_\ulamb({\cal E})$:
therefore the
${\cal E}$-Schubert varieties $\overline{W_\ulamb({\cal
E})}$ form a {\em cellular decomposition} of $G$ in the sense
of~\cite{Fu1}, p. 23, Example 1.9.1. By the same reference, their classes
modulo linear equivalence, $\ovW_\ulamb$, freely generate
$A_*(G)$ over the integers. This explain formula~(\ref{eq:shc1}) used to define
the Schubert (co)cycles $\sigma_\ulamb$.

If $h\in\NN$, set $\sigma_h=(h,0,\ldots,0)$, then
$\sigma_\ulamb=\sigma_h$: it is said to be a {\em special Schubert (co)cycle}
(\cite{GH}, p.~203). The presentation of the Chow ring of $G$ is (see
e.g.~\cite{Bottu}, p.~293):
\be
A^*(G)={\ZZ[\sigma_1,\ldots,\sigma_{n-k}]\over
(Y_{k+1}(\sigma),\ldots,Y_n(\sigma))}.\label{eq:prescrg}
\ee
where $Y_i(\sigma)$'s are defined by the formal equality:
\be
{1\over 1+\sigma_1t+\sigma_2t^2+\ldots+\sigma_{n-k}t^k}=\sum_{i\geq
0}(-1)^iY_i(\sigma)t^i.\label{eq:formequal}
\ee
In particular, $A^*(G)$ is generated by
the  {\em special Schubert cycles} only. The aim of Schubert Calculus is to
compute the structure constants $\{C^\unu_{\ulamb\umu}\}$ of such an
algebra:
\be
\sigma_\ulamb\cup\sigma_\umu=\sum_{|\unu|=|\ulamb|+|\umu|}
C^\unu_{\ulamb\umu}\sigma_\unu.
\ee
The structure constants $\{C^\unu_{\ulamb\umu}\}$ can be algorithmically
determined once one knew all the products $\sigma_h\cup\sigma_\ulamb$ and the
shape of all the $\sigma_\ulamb$'s as polynomial expressions in the special
Schubert cycles. These are prescribed by the well known equalities in $A^*(G)$:
\begin{enumerate}
\item Pieri's formula (Cf. \cite{Fu1}, p.~271 ;\cite{GH}, p.~203):
\be
\sigma_h\cup\sigma_\ulamb=\sum_{\umu}\sigma_\umu,\label{eq:piericlas}
\ee
the sum extended over all $\umu$ such that $\mu_1\geq \lambda_1\geq
\mu_2\geq\lambda_2\geq\ldots\geq \mu_k\geq\lambda_k$ and
$|\umu|=|\ulamb|+h$,
and
\item  Giambelli's formula (Cf. \cite{Fu1}, p.~271; \cite{GH}, p.~205):
\be
\sigma_\ulamb=\left|\matrix{\sigma_{\lambda_k}&
\sigma_{\lambda_{k-1}+1}&\ldots&\sigma_{\lambda_{1}+k-1}\cr
\sigma_{\lambda_k-1}&\sigma_{\lambda_{k-1}}&\ldots&\sigma_{\lambda_{1}+k-2}\cr
\vdots&\vdots&\ddots&\vdots\cr\sigma_{\lambda_k-k+1}
&\sigma_{\lambda_{k-1}-k+2}&\ldots&\sigma_{\lambda_1}}\right|=
\Delta_\ulamb(\sigma).
\ee
the product being the {\em cup product}.
\end{enumerate}
These formulas also determine the  structure of $A_*(G)$ as a module over
$A^*(G)$:
\begin{enumerate}
\item
\be
(\sigma_h\cup\sigma_\ulamb)\cap [G]=\sigma_i\cap(\sigma_\ulamb\cap [G])=
\sigma_h\cap \overline{W}_\ulamb=\sum_{\umu}\overline{W}_\umu
\ee
the sum being over all $\umu$ such that $\mu_1\geq \lambda_1\geq
\mu_2\geq\lambda_2\geq\ldots\geq \mu_k\geq\lambda_k$ with
$|\umu|=|\ulamb|+h$,
and
\item
\be
\overline{W}_\ulamb=\Delta_\ulamb(\sigma)\cap [G].
\ee

\end{enumerate}
In \cite{Wi}, Witten puts a {\em quantum product} $\cup_q$
on
$A^*(G)\otimes_{\bf Z}\ZZ[q]$, making it into a commutative graded ring
with unit (the fundamental class $[G]$ of $G$), here shortly denoted 
as $QA^*(G)$. Such a {\em quantum cup} product is defined as:
\[
\sigma_\ulamb\cup_q\sigma_\umu=\sum_{d\geq 0}\,\,\sum_{|\unu|=
|\ulamb|+|\umu|-nd}C^\unu_{\ulamb\umu}(d)q^d\sigma_\unu,
\]
where $C^\unu_{\ulamb\umu}(d)$ is the number of maps
$f:\PP^1\lra G$ of degree $d$ such that $f(0)\in\overline{W_{\ulamb}({\cal
E})}$, $f(1)\in\overline{W_{\umu}({\cal
F})}$ and $f(\infty)\in\overline{W_{\unu^\vee}({\cal
G})}$, where ${\cal E}$, ${\cal F}$, ${\cal G}$ are three flags in sufficiently
general position. Let $QA_*(G):=A_*(G)\otimes_{\bf Z}\ZZ[q]$ be the
$\ZZ[q]$-module generated by the Schubert cycles $\overline{W}_\umu$. Then
$QA_*(G)$ is a module over $QA^*(G)$ via  {\em quantum cap product} $\cap_q$,
defined on generators by
\[
\left\{\matrix{\sigma_\ulamb\cap_q [G]&=&\sigma_\ulamb\cap [G]=\ovW_\ulamb  \cr
\sigma_\ulamb\cap_q
\ovW_\umu&=&(\sigma_\ulamb\cup_q \sigma_\umu)\cap [G]}\right. ,
\]
and extended by $\ZZ[q]$-linearity. The presentation of the
Quantum intersection ring $QA^*(G)$ has been computed by Witten and
Siebert-Tian (\cite{Wi},\cite{SiebTi}):
\be
QA^*(G):={\ZZ[q][\sigma_1,\ldots,\sigma_{n-k}]\over{(Y_{k+1}(\sigma),\ldots,
Y_n(\sigma)-(-1)^{n-k-1}q})},\label{eq:wstring}
\ee
where still $Y_i(\sigma)$ are defined by the formal
equality~(\ref{eq:formequal}) but computed according to the quantum product.
Therefore the structure constants
$\{C^\unu_{\ulamb\umu}(d)\}$ of $QA^*(G)$ are determined once one knows some
quantum analogue of Pieri's and Giambelli's formula. Such formulas, listed
below, have been found and proven by Bertram in~\cite{Ber1} (see
also~\cite{Buch} for a simpler alternative proof):
\begin{enumerate}
\item{\em Quantum
Pieri's formula}:
\be
\sigma_h\cup_q\sigma_\ulamb=\sigma_h\cup\sigma_\ulamb+q\sum_{\umu}\sigma_\umu,
\ee
the sum over all $\umu$ such that $\mu_1\geq\lambda_1-1\geq
\mu_2\geq\lambda_2-1\geq\ldots\geq\mu_k\geq \lambda_k-1\geq 0$,
$|\mu|=|\ulamb|+h-n$,
and
\item {\em Quantum Giambelli's formula}:
\be
\sigma_\ulamb=\Delta_{\ulamb}(\sigma),
\ee
\end{enumerate}
where Giambelli's determinant is now computed using the quantum cup product
$\cup_q$.

\section{The Algebra of Schubert Calculus}\label{Sect3}

\subsection{Hasse-Schmidt Derivations on Exterior Algebras}\label{sect1.2}

\begin{defin}
Let $\wM$ be the exterior algebra of a module over an integral domain
$A$. A {\em
Hasse-Schmidt} ($HS$-)derivation $D_t$ on $\wM$ is an $A$-algebra homomorphism:
\[
D_t:\wM\lra\left(\wM\right)[[t]].
\]
\end{defin}

A Hasse-Schmidt derivation on $\wM$ determines, and is determined by,
its {\em coefficients} $D_i:\wM\lra\wM$, defined via the equality
\[
D_t\alpha=\sum_{i\geq 0}D_i(\alpha)t^i,\qquad \forall\alpha\in\wM.
\]
\begin{prop} The endomorphisms $D_i$ satisfy  (generalized) Leibniz's rule:
\be
D_i(\alpha\wedge\beta)=\sum_{\matrix{i_1+i_2=i\cr i_1\geq0,
i_2\geq 0}}D_{i_1}\alpha\wedge D_{i_2}\beta,\label{eq:for2}
\ee
\end{prop}
\proof
Since $D_t$ is an $A$-algebra homomorphism, one has:
\[
D_t(\alpha\wedge\beta)=D_t\alpha\wedge D_t\beta,\quad
\forall\alpha,\beta\in\wM.
\]
Therefore $D_i(\alpha\wedge\beta)$ is the coefficient of $t^i$ in the expansion
of $D_t(\alpha\wedge\beta)$, which is also the coefficient of $t^i$ in the
expansion of the wedge product
\[
(D_0\alpha+ D_1\alpha\, t+ D_2\alpha\, t^2+\ldots)\wedge (D_0\beta+
D_1\beta \,t+ D_2\beta\, t^2+\ldots),
\]
i.e. exactly the right hand side of eq.~(\ref{eq:for2}).

\qed

In particular, $D_0$ is an algebra homomorphism and $D_1$ a usual derivation:
\[
D_0(\alpha\wedge\beta)=D_0(\alpha)\wedge D_0(\beta),\qquad
D_1(\alpha\wedge\beta)=D_1\alpha\wedge\beta+\alpha\wedge D_1\beta.
\]

If $D_t=\sum_{i\geq 0}D_it^i$, the sequence of the coefficients
\[
D:=(D_0,D_1,D_2,\ldots)
\]
of $D_t$ will be  also said
a $HS$-derivation and the $D_i$'s will be also said to be the {\em
components} of
$D$. To denote a $HS$-derivation the symbols $D$ and $D_t$ shall be used
interchangeably.

\begin{defin}
A $HS$-derivation $D$ on $\wM$ is said to be {\em regular} if $D_0\in
End_A(\wM)$ is an $A$-automorphism; a regular $HS$-derivation is {\em
normalized} if
$D_0=id_{\bigwedge M}$.
\end{defin}

\begin{prop}
Let $HS_t(\wM)$ be the set of all regular $HS$-derivations of $\wM$. Then
$HS_t(\wM)$ is a group.
\end{prop}
\proof First of all, if $D_t,D'_t\in HS_t(\wM)$, one has:
\begin{eqnarray*}
D_t\circ
D_t'(\alpha\wedge\beta)&=&D_t(D_t'(\alpha\wedge\beta))=D_t(D_t'(\alpha)\wedge
D'_t(\beta))=\\&=&D_t(D_t'(\alpha))\wedge
D_t(D'_t(\beta))=D_t\circ
D_t'(\alpha)\wedge D_t\circ
D_t'(\beta),
\end{eqnarray*}
so that $D_t\circ D_t'\in HS_t(\wM)$. The composition $\circ$ is obviously
associative. The identity $1:=1_M$ of $M$ belongs to $HS_t(\wM)$.  Let $E_t$
be the formal inverse of $D_t$, (i.e. $D_t\circ E_t=E_t\circ 
D_t=id_{\wM}$), existing by
the regularity hypothesis on $D_t$. One has:
\begin{eqnarray*}
E_t(\alpha\wedge\beta)&=&E_t(D_tE_t(\alpha)\wedge 
D_tE_t(\beta))=E_t(D_t(E_t\alpha\wedge 
E_t\beta))=\\
&=&(E_tD_t)(E_t\alpha\wedge
E_t\beta)=(E_t\alpha\wedge
E_t\beta),
\end{eqnarray*}
i.e. $E_t\in HS_t(\wM)$, and is obviously regular.
\qed

\medskip
\noindent
If $A\lra B$ is a ring homomorphism and $D_t$ is a $HS$ derivation on
$\wM$, the $B$-linear extension of $D_t$ to $\wM\otimes_AB$ is
  a
$HS$-derivation, too, denoted by the same
symbol
$D_t$, abusing notation.

\subsection{Intersection Pairs} \label{sect1.4}
Let $M$ be a free module over a n\"otherian integral domain $A$. Let $D$
be a sequence $(D_0,D_1,\ldots)$ of pairwise commuting $A$-endomorphisms of $M$
($D_iD_j=D_jD_i$), said to be the {\em components} of $D$. It may also be
represented as a homomorphism $D_t\in End_A(M)[[t]]$, putting $D_t=\sum_{i\geq
0}D_it^i$: in this case we shall refer to the $D_i$'s as to the {\em 
components}
of
$D_t$.
Let $A[\TT]:=A[T_1,T_2,\ldots]$ be the $A$-polynomial ring in infinitely many
indeterminates, graded by declaring that $\deg(T_i)=i$, and let $A[D]$ be the
commutative subring of $End_A(M)$ generated by the $D_i$.  Clearly
$M$ is an
$A[D]$-module. Let $E_D:A[\TT]\lra
A[D]$ be the natural evaluation epimorphism, defined by
$E_D(P(\TT))=P(D)$ (``evaluating $T_i$ at $D_i$") for each $P(\TT)\in A[\TT]$.
If $m\in M$, let $E_m:A[\TT]\lra M$ be the evaluation
morphism $E_m(P(\TT))=E_D(P(\TT))\cdot m=P(D)\cdot m$.
\begin{defin} The data $(M,D_t)$ is said to be an $A$-{\em intersection
pair} ($A-{\cal IP}$) if there is $m_0\in M$, said to be a {\em fundamental
element}, such that
$E_{m_0}$ is an epimorphism. A $\ZZ-{\cal IP}$ will be shortly said {\em
intersection pair}.
\end{defin}

Obviously $\ker E_{m_0}=\ker E_D$. Therefore $A[\TT]/\ker(E_{m_0})\cong
A[D]\cong M$. The ring ${\cal A}^*(M,D_t):=A[D]$ is said to be the
{\em intersection ring} of the intersection pair $(M,D_t)$. The isomorphism:
\[
\matrix{{\Pi}:&{\cal A}^*(M,D_t)&\lra&M\cr{ }&P(D)&\longmapsto&P(D)\cdot m_0}
\]
will be said to be the {\em Poincar\'e isomorphism}, while the ring
$A[\TT]/\ker(E_{m_0})$ is the {\em presentation} of $A[D]$.
Let ${\cal A}^h(M,D_t)=A[D]_h=E_D(A[\TT]_h)$, where
\be
A[\TT]_h=\bigoplus_{\matrix{m_1+2m_2\ldots+hm_h=h\cr m_i\geq 0}} A\cdot
\TT_1^{m_1}\cdot\ldots\cdot \TT_h^{m_h}=\bigoplus_{|\ulamb|=h} A\cdot
T_\ulamb,\label{eq:forgrading}
\ee
and $T_\ulamb=T_{1}^{m_1}\cdot\ldots\cdot
T_{h}^{m_h}$, having set $\ulamb=(1^{m_1}\ldots h^{m_h})$. Then $A[D]$
induces a grading on $M$:
\[
M=M_0\oplus M_1\oplus M_2\oplus\ldots
\]
by setting $M_h=A[D]_h\cdot m_0$. A homogeneous element $m\in M_h$ is 
said to be
of {\em weight} $h$.

\begin{prop}
Let $(M,D_t)$ be an $A-{\cal IP}$ and let $A'$ be an $A$-algebra. Let $D'$ be
the
$A'$-linear extension of $D$ to $M'=M\otimes_AA'$. Then $(M',D')$ is 
a $A'-{\cal
IP}$, having the same fundamental element as $(M,D_t)$. Furthermore ${\cal
A}^*(M',D_t')={\cal A}^*(M,D_t)\otimes_AA'$.
\end{prop}
\proof

First of all $M'$ is clearly a free  $A'$-module.
Any $m'\in M'$ is a finite linear $A'$-combination
$a'_1m_1+\ldots+a'_hm_h\in M'$. Then there exists $G_{m_i}\in A[\TT]$ such that
$m_i=G_{m_i}(D)m_0$, where
$m_0$ is the fundamental element of the $A-{\cal IP}$ $(M, D_t)$. Then
$G_{m'}:=a'_1G_{m_1}+\ldots+a'_hG_{m_h}\in A'[\TT]$ is a polynomial such that
$E_D(G_{m'})\cdot m_0=m'$. This proves that $(M',D')$ is an $A'-{\cal IP}$
with fundamental element $m_0\in M$ -- the fundamental element of $(M, D_t)$.
Clearly ${\cal A}^*(M,D_t)\otimes_AA'\cong M\otimes_AA'\cong M'$, so 
that ${\cal
A}^*(M',D')={\cal A}^*(M,D_t)\otimes_AA'$.

\qed

Let $\bt=(t_1,\ldots, t_k)$ be a set of $k$ indeterminates. Then, for each
$1\leq i\leq k-1$
$(M[[t_1,\ldots t_{i-1}]], D_{t_i})$ is an $A[[t_1,\ldots t_{i-1}]]-{\cal IP}$,
once $D_{t_i}=\sum_{j\geq 0}D_jt_i^j$ has been extended to $M[[t_1,\ldots
t_{i-1}]],$  by $A[[t_1,\ldots t_{i-1}]]$-linearity.

Therefore it makes sense to consider the
composition
\[
D_{t_1}\circ\ldots\circ D_{t_k}:M\lra M[[t_1,\ldots,t_k]]
\]

\[
\matrix{M&\lra& M[[t_k]]&\lra &M[[t_{k-1},t_k]]&\lra&\ldots&\lra&
M[[t_1,\ldots,t_k]]\cr
{ }&D_{t_k}&{ }&D_{t_{k-1}}&{ }&D_{t_{k-2}}&\ldots&{D_{t_1}}&{ } }
\]
Then:
\begin{lemma}\label{lemmapreGia}
The following formula holds:
\be
D_{t_1}\circ\ldots\circ D_{t_k}=\sum_{\ulamb\in{\cal L}_{\leq k}}
s_\ulamb(T)\Delta_\ulamb(D),\label{eq:HSS61}
\ee
where $s_\ulamb(T)$ are the Schur symmetric polynomials defined in
Sect.~\ref{sect1.2}.
\end{lemma}

\proof

Make the substitution
$D_i=h_i(X)$, where $h_i(X)$ is the complete symmetric polynomial of
degree
$i$ in
$k$ formal variables
$X=(x_1,\ldots, x_k)$. Then
$D_{t_i}$ can be expressed as:
\[
D_{t_i}=\prod_{j=1}^k{1\over 1-x_jt_i}.
\]
It follows that:
\be
D_{t_1}\circ\ldots\circ D_{t_k}=\prod_{i,j=1}^k{1\over
1-x_jt_i}=\sum_{\ulamb\in{\cal L}_{\leq k}} s_\ulamb(T)s_\ulamb(X),
\label{eq:HSS6}
\ee
where last equality of formula (\ref{eq:HSS6}) follows from
identity (4.3) at p. 33 of \cite{MacDon}. But in the old
variables,
$s_\ulamb(X)$ is exactly the determinantal formula relating the
$s_\ulamb$ with the complete symmetric polynomials (see formula~(3.4) of
\cite{MacDon}, p.~25), which are the
$D_i's$, i.e. exactly formula~(\ref{eq:HSS61}).
\qed

\subsection{The Canonical Intersection Pair }\label{Sect3.1}
Let $(M,D_t)$ be the data of
\begin{enumerate}
\item \label{CP1} an infinite free $\ZZ$-module $M=\oplus_{i\geq
1}\ZZ\epsilon^i$
   generated by $(\epsilon^1,\epsilon^2,\ldots)$ (each element
of $M$ is a finite free $\ZZ$-linear combination of the $\epsilon^i$'s);
\item \label{CP2}
a formal power series $D_t:=\sum_{i\geq 0}D_it^i\in End_{\bf Z}(M)[[t]]$ whose
{\em coefficients}
  are the {\em components} of
a  sequence of endomorphisms
\be
D=(D_0,D_1,D_2,\ldots)\label{eq:comps}
\ee
such
that $D_0$ is the identity, $D_1$ is the one-step shift operator, i.e.
$D_1(\epsilon^i)=\epsilon^{i+1}$, and $D_i=D_1^i$.
\end{enumerate}

\begin{prop}\label{propCP}
The pair $(M,D_t)$ is an ${\cal IP}$ said to be {\em canonical intersection
pair},  and its intersection ring is isomorphic to $\ZZ[D_1]$.
\end{prop}
\proof
Let $M_i=\ZZ\cdot\epsilon^i$. Then $M=\oplus M_i$ and $D_j:M\lra M$ is a
homogeneous endomorphism of degree  $j$ ($D_jM_i\subset M_{i+j}$).

If $E_D:\ZZ[\TT]\lra End_{\bf Z}(M)$  is the natural
evaluation homomorphism gotten by evaluating $\phi(\TT)\in\ZZ[\TT]$ at
$T_i=D_i$ (Sect.~\ref{sect1.4}), the image $\ZZ[D]:=Im(E_D)$ is a
commutative subalgebra of
$End_{\bf Z}(M)$, because
$D_iD_j=D_jD_i$. \linebreak Let $E_{\epsilon^1}:\ZZ[\TT]\lra M$ be the
evaluation morphism sending $\phi(\TT)\in\ZZ[\TT]$ onto $\phi(D)\epsilon^1$.
Then
\[
\left\{\matrix{&{E}_{\epsilon^1}:&\ZZ[\TT]&\lra&
M\cr
&{}&\phi(\TT)&\longmapsto&\phi(D)\epsilon^1}\right.
\]
is obviously surjective. In fact, for any
$m=a_{i_1}\epsilon^{i_1}+\ldots+a_{i_j}\epsilon^{i_j}$, one has:
\[
{E}_{\epsilon^1}(a_{i_1}T_{i_1-1}+\ldots+a_{i_j}T_{i_j-1})=
(a_{i_1}D_{i_1-1}+\ldots+a_{i_j}D_{i_j-1})\epsilon^1=m.
\]
Therefore $M$ is isomorphic, as a $\ZZ$-module, to the ring
\[
{\cal
A}^*(M,D_t):={\ZZ[T]\over\ker(\tilde{E}_{\epsilon^1})}.
\]
One easily checks that
$\ker(\tilde{E}_{\epsilon^1})=(T_i-T_1^i)_{i\geq 1}$. Hence:
\[
{\cal
A}^*(M,D_t)={\ZZ[T_1,T_2\ldots]\over (T_2-T_1^2,T_3-T_1^3,\ldots
)}\cong\ZZ[T_1]\cong E_D(\ZZ[T_1])=\ZZ[D_1].
\]

\qed

\subsection{Schubert Derivations}
Let $k\geq 1$. Given  a {\em canonical intersection pair} $(M,D_t)$
(Sect.~\ref{Sect3.1}) let
$\wM=\oplus_{k\geq 0}\bigwedge^kM$  be the exterior algebra of the
$\ZZ$-module
$M$. Denote by the same symbol the extension of $D_t$ to all of $End_{\bf
Z}(\wM)[[t]]$. In such a way $D_t$ turns into  a
$HS$-derivation on $\wM$, i.e.:
\be
D_t(\alpha\wedge\beta)=D_t\alpha\wedge D_t\beta,\quad \forall
\alpha,\beta\in \wM.\label{eq:for1}
\ee
\begin{defin}
The formal power series $D_t$ so obtained will be said {\em Schubert
derivation} or, briefly, ${\cal S}$-derivation.
\end{defin}
  A ${\cal S}$-derivation
is a normalized $HS$-derivation. In fact, for each $\ikformep$,
\[
D_0(\ikformep)=D_0\epsilon^{i_1}\wedge\ldots\wedge D_0\epsilon^{i_k}=\ikformep,
\]
because $D_0$ is the identity of $M$.
Moreover, for $i\geq 1$, each $D_i$ satisfies the generalized 
Leibniz's rule, as
explained in Sect.~\ref{sect1.2}. The sequence
$D:=(D_0,D_1,D_2,\ldots)$ of the coefficients of $D_t$ will also be called
${\cal S}$-derivation.
\begin{prop}
The endomorphisms $D_i:\wkM\lra \wkM$ are pairwise commuting.
\end{prop}
\proof
By induction on $k$. For $k=1$ the claim is true by construction. Assume that
the property holds for $k-1$. Since any $m\in \wkM$ is a finite sum of $k$
vectors of the form $\alpha\wedge \beta$, with $\alpha\in M$ and $\beta\in
\bigwedge^{k-1}M$, without loss of generality one may check the property for
any
$m$ of this form.
One has:
\be
D_iD_j(\alpha\wedge\beta)=D_i\left(\sum_{j_1+j_2=j}D_{j_1}\alpha\wedge
D_{j_2}\beta\right)=\sum_{i_1+i_2=i}\,\,\sum_{j_1+j_2=j}D_{i_1}D_{j_1}\alpha\wedge
D_{i_2}D_{j_2}\beta.\label{eq:for8}
\ee
By the inductive hypothesis, last member of equality~(\ref{eq:for8}) is equal
to:
\[
\sum_{j_1+j_2=j}\,\,\sum_{i_1+i_2=i}D_{j_1}D_{i_1}\alpha\wedge
D_{j_2}D_{i_2}\beta=D_j\left(\sum_{i_1+i_2=i}D_{i_1}\alpha\wedge
D_{i_2}\beta\right)=D_jD_i(\alpha\wedge\beta).
\]
\qed

As a consequence the natural evaluation morphism $E_D:\ZZ[\TT]\lra End_{\bf
Z}(\wkM)$ maps $\ZZ[\TT]$ onto a commutative subring $\ZZ[D]$ of $End_{\bf
Z}(\wkM)$.

As for more general $HS$-derivations on $\wM$:
\begin{prop}\label{propprepie}
The ``generalized" Leibniz's rule holds:
\[
D_h(\epsilon^{i_1}\wedge \epsilon^{i_2}\wedge\ldots\wedge \epsilon^{i_k})=
\sum_{\matrix{h_1+\ldots+h_k=h\cr h_i\geq 0}} \epsilon^{i_1+h_1}\wedge
\epsilon^{i_2+h_2}
\wedge\ldots\wedge \epsilon^{i_k+h_k}.
\]
\end{prop}
\proof
By induction on $k$. If $k=1$ the property is obviously true. Now
assume it holds for $k-1$. Then
\be
D_h( \epsilon^{i_1}\wedge \epsilon^{i_2}\wedge\ldots\wedge \epsilon^{i_k})=
\sum_{h_1=0}^{h} \epsilon^{i_1+h_1}\wedge
D_{h-h_1}( \epsilon^{i_2}\wedge\ldots\wedge \epsilon^{i_k}).
\label{eq:prepier}
\ee
By the inductive hypothesis:
\[
D_{h-h_1}( \epsilon^{i_2}\wedge\ldots\wedge \epsilon^{i_k})=
\sum_{h_2+\ldots+h_k=h-h_1} \epsilon^{i_2+h_2}\wedge\ldots\wedge
  \epsilon^{i_k+h_k}.
\]
Hence the right hand side of formula~(\ref{eq:prepier}) turns into:
\[
D_h( \epsilon^{i_1}\wedge \epsilon^{i_2}\wedge\ldots\wedge \epsilon^{i_k})=
\sum_{h_1+\ldots+h_k=h} \epsilon^{i_1+h_1}\wedge\ldots
\wedge \epsilon^{i_k+h_k}.
\]
\qed

\begin{example}\label{expier}{\rm Let us
compute
$D_2(\epsilon^2\wedge \epsilon^3\wedge \epsilon^5)\in \bigwedge^3M$. One has:
\begin{eqnarray*}
&{ }&D_2(\epsilon^2\wedge \epsilon^3\wedge \epsilon^5)=D_2\epsilon^2\wedge
\epsilon^3\wedge \epsilon^5+ D_1\epsilon^2\wedge D_1(\epsilon^3\wedge
\epsilon^5)+ \epsilon^2\wedge D_2(\epsilon^3\wedge \epsilon^5)=\\
\\&{}& =\epsilon^4\wedge \epsilon^3\wedge \epsilon^5 + \epsilon^3\wedge
(\epsilon^4\wedge \epsilon^5+\epsilon^3\wedge \epsilon^6)+ \epsilon^2\wedge
(\epsilon^5\wedge \epsilon^5+ \epsilon^4\wedge \epsilon^6+\epsilon^3\wedge
\epsilon^7)=\\ \\ &{}&=\epsilon^4\wedge \epsilon^3\wedge \epsilon^5 +
\epsilon^3\wedge
\epsilon^4\wedge \epsilon^5+\epsilon^3\wedge \epsilon^3\wedge \epsilon^6+
\epsilon^2\wedge \epsilon^5\wedge \epsilon^5 + \epsilon^2\wedge 
\epsilon^4\wedge
\epsilon^6+\\&{}&\hskip 70pt +\epsilon^2\wedge \epsilon^3\wedge \epsilon^7=
\epsilon^2\wedge \epsilon^4\wedge \epsilon^6+\epsilon^2\wedge \epsilon^3\wedge
\epsilon^7,
\end{eqnarray*}
where last equality results from the vanishing of the terms
having equal factors and to the cancelations due to skew-symmetry.
}
\end{example}
Example~\ref{expier} shows that computing ${\cal S}$-derivatives of
$k$-vectors is a straightforward matter and that one has to care only
about possible vanishing and cancelations.

However, the practise of many examples naturally suggests the following:

\begin{prop}\label{Pieri} {\em Pieri's formula for ${\cal
S}$-derivatives} holds:
\[
D_h(\ikformep)=\sum \epsilon^{j_1}\wedge\ldots\wedge \epsilon^{j_k},
\]
where the sum is over all $i_1\leq j_1< i_2\leq j_2<\ldots<
i_k\leq j_k$, such that $\sum_{p=1}^k
j_p=h+\sum_{p=1}^k i_p$.

\end{prop}
\proof Let $I=(1\leq i_1<\ldots<i_k)$. By 
Proposition~\ref{propprepie}, one has:
\be
D_h(\ikformep)=\sum_{\matrix{h_1+\ldots+h_k=h\cr h_i\geq
0}}\epsilon^{i_1+h_1}\wedge\ldots\wedge \epsilon^{i_k+h_k}.\label{eq:sumpie}
\ee
Now, each term such that $i_p+h_p=i_q+h_q$, for some $1\leq p<q\leq k$,
vanishes, by  skew-symmetry of the wedge product. Moreover, if in the
expansion~(\ref{eq:sumpie})  a term like
\be
\epsilon^{i_1+h_1}\wedge\ldots\wedge \epsilon^{i_p+h_p}\wedge\ldots\wedge
  \epsilon^{i_q+h_q}\wedge\ldots\wedge
\epsilon^{i_k+h_k},\label{eq:1term}
\ee
with
\be
i_p+h_p>i_q+h_q,\label{eq:ineqpie}
\ee
occurs,
then, keeping the same
$h_j$ for $j\neq p,q$ in the other factors, it will also occur a term
of the form
\be
\epsilon^{i_1+h_1}\wedge\ldots\wedge \epsilon^{i_p+h'_p}\wedge\ldots\wedge
  \epsilon^{i_q+h'_q}\wedge\ldots\wedge
\epsilon^{i_k+h_k},\label{eq:-1term}
\ee
where $h'_p=(i_q-i_p)+h_q$  and
$h'_q=(i_p-i_q)+h_p$.  In fact  $h'_p$ is
obviously non negative, $h'_q> 0$ by virtue
of inequality~(\ref{eq:ineqpie}) and
$h'_p+h'_q+\sum_{j\neq p,q}h_j=h_p+h_q+\sum_{j\neq p,q}h_j=h$. Since
$i_p+h'_p=i_q+h_q$ and
$i_q+h'_q=i_p+h_p$, the terms~(\ref{eq:1term}) and~(\ref{eq:-1term})
occurring in the sum~(\ref{eq:sumpie})  cancel out, proving the claim.
Therefore, in the sum~(\ref{eq:sumpie}), only the summands with
\[
i_1\leq i_1+h_1<i_2\leq i_2+h_2<\ldots <i_{k-1}\leq
i_{k-1}+h_{k-1}< i_k\leq i_k+h_k
\]
  survive.
Putting $J=(j_1,\ldots,j_k)=(i_1+h_1,\ldots,i_k+h_k)$,
one finally gets:
\be
D_h(\ikformep)=\sum\epsilon^{j_1}\wedge\ldots\wedge
\epsilon^{j_k},\label{eq:pier1}
\ee
where the sum is over all $J$ such that $i_1\leq j_1<i_2\leq j_2<\ldots
<i_k\leq j_k\leq i_k+h$ and $|J|=|I|+h$
\qed

\medskip
As an application of Pieri's formula, one can observe that:
\begin{corol}~\label{propN2}
The following equality holds in $\bigwedge^kM$:
\be
D_{h}(\epsilon^{i}\wedge \epsilon^{i+1}\wedge\ldots\wedge
  \epsilon^{i-2+k}\wedge \epsilon^{i-1+k})= \epsilon^{i}\wedge
\epsilon^{i+1}\wedge\ldots\wedge \epsilon^{i-2+k}\wedge
  \epsilon^{i-1+k+h}.
\ee
\end{corol}
\proof

Apply Proposition~\ref{Pieri}:
\[
D_{h}(\epsilon^{i}\wedge \epsilon^{i+1}\wedge\ldots\wedge
  \epsilon^{i+k-2}\wedge \epsilon^{i+k-1})=\sum \epsilon^{j_1}\wedge\ldots\wedge
\epsilon^{j_k}
\]
with  $|J|=|I|+h$ and $i\leq j_1<i+1\leq j_2<i+2\leq \ldots <i+k-2\leq
j_{k-1}<i+k-1\leq j_k$. Hence $j_p=i+p-1$, for all $1\leq p\leq k-1$,
  and
$j_k=i-1+k+h$, as claimed.
\qed

\begin{thm}\label{giamth} Let $\ulamb=(r_kr_{k-1}\ldots r_1)\in {\cal L}_{\leq
k}$. The morphism
\[
E_{\kformep}:\ZZ[D]\lra \bigwedge^kM
\]
is an epimorphism. In fact,
{\em Giambelli's formula} holds:
\be
\Delta_\ulamb(D)\kformep=\rikformep\label{eq:gengiamb}
\ee
\end{thm}
\proof
Let $\bt:=(t_1,\ldots, t_k)$ be a set of $k$ indeterminates over $\ZZ$.
Then, by composition, one may define (see Section~\ref{sect1.4}):
\[
D_{t_1}\circ\ldots\circ
D_{t_k}:\bigwedge^kM\lra\bigwedge^kM[[t_1,\ldots,t_k]].
\]
Clearly:
\be
(D_{t_1}\circ\ldots\circ
D_{t_k})(\kformep)=(D_{t_1}\circ\ldots\circ
D_{t_k})(\epsilon^{1})\wedge\ldots\wedge (D_{t_1}\circ\ldots\circ
D_{t_k})(\epsilon^{k})\label{eq:HSS3}
\ee

Applying~\ref{lemmapreGia} to both sides of~(\ref{eq:HSS3}), one has:

\be
\sum_{\ulamb\in{\cal
L}_{\leq
k}}s_\ulamb(\bt)\Delta_\ulamb(D)\kformep=
\sum_{\ulamb\in{\cal
L}_{\leq
k}}s_\ulamb(\bt)\Delta_\ulamb(D)\epsilon^{1}\wedge\ldots\wedge
\sum_{\ulamb\in{\cal
L}_{\leq
k}}s_\ulamb(\bt)\Delta_\ulamb(D)\epsilon^{k}
\ee
The right hand side is equal to

\begin{eqnarray}
&{ }&\sum_{r_1\geq
0}s_{r_1}(\bt)\Delta_{r_1}(D)\epsilon^1\wedge\ldots\wedge
\sum_{r_k\geq
0}s_{r_k}(\bt)\Delta_{r_k}(D)\epsilon^k\nonumber\\&{ }&=\sum_{r_1\geq
0}h_{r_1}(\bt)\epsilon^{1+r_1}\wedge\ldots\wedge
\sum_{r_k\geq 0}h_{r_k}(\bt)\epsilon^{k+r_k}=\nonumber\\&{}&=\sum_{\ulamb\in
{\cal L}_{\leq k}} s_\ulamb(\bt)\cdot
  \epsilon^{1+r_1(\ulamb)}\wedge\ldots\wedge
\epsilon^{k+r_k(\ulamb)},\label{eq:HSS4}
\end{eqnarray}
last equality coming from a straightforward computation of an
exterior product and the determinantal relationship between the
complete symmetric polynomials $h_i(\bt)$ and the Schur symmetric
polynomials
$s_\ulamb(\bt)$ (Sect. \ref{sect2.2}).

  Since $(s_\ulamb(\bt):\ulamb\in {\cal
L}_{\leq k})$ is a basis of the symmetric polynomials in the
indeterminates $\bt$ with integral coefficients, formula~(\ref{eq:gengiamb})
follows.
\qed

\subsection{The intersection ring of $(\wkM,D)$}

\subsubsection{Generalities}
As usual, denote by $\TT$ the set of infinitely many indeterminates $(T_1, T_2,
T_3,\ldots)$.

\begin{corol}
The polynomials $\{\Delta_\ulamb(\TT), \ulamb\in {\cal P}\}$ are
$\ZZ$-linearly independent.
\end{corol}
\proof
Suppose that
$
\sum_{\ulamb\in {\cal P}'}a_\ulamb\Delta_\ulamb(\TT)=0
$
is a non trivial linear dependence relation,
where ${\cal P}'$ is a finite subset of ${\cal P}$.
Then there is $k\geq 1$, such that ${\cal P}'\subseteq {\cal L}_{\leq k}$ and
\[
\sum_{\ulamb\in {\cal P}'}a_\ulamb\Delta_\ulamb(D)\kformep=0,
\]
which is impossible, because
$\Delta_\ulamb(D)\kformep=\epsilon^{1+r_1(\ulamb)}\wedge\ldots\wedge
\epsilon^{k+r_k(\ulamb)}$
are all linearly independent, being part of a basis of $\bigwedge^kM$.
\qed
\begin{lemma}\label{deltabasis}
The subset $\Delta^h(\TT)=\{\Delta_\ulamb(\TT):\ulamb\in {\cal P}_h\}\subset
\ZZ[\TT]$, is a
$\ZZ$-basis of $\ZZ[\TT]_h$.
\end{lemma}
\proof
For any $\ulamb\in {\cal P}_h$, $\Delta_\ulamb(\TT)$ is a unique $\ZZ$-linear
combination of $T_\umu$'s, $(\umu\in {\cal P}_h)$:
\[
\Delta_\ulamb(\TT)=\sum_{\umu_{\ulamb}\in{\cal
P}_h}a_{\umu_{\ulamb}}D_{\umu_{\ulamb}}.
\]
Then, since $\sharp\{\Delta^h(\TT)\}=\sharp\{T_\umu:\umu\in{\cal P}_h\}$, and
all the $\Delta_\ulamb\in\Delta^h(D)$ are linearly independent, it follows that
they freely span $\QQ[\TT]_h=\ZZ[\TT]_h\otimes_{\bf Z}\QQ$. In particular, for
each $\umu\in{\cal P}_k$,
$T_\umu$ is a unique $\QQ$-linear combination of
$\Delta_\ulamb(\TT)\in\Delta^h(\TT)$:
\be
T_\umu=\sum_{\ulamb_{\umu}}b_{\ulamb_{\umu}}\Delta_{\ulamb_{\umu}}(\TT).\label{eq:for5}
\ee
We only need to show that the coefficients $b_{\ulamb_{\umu}}\in\ZZ$. But
equality~(\ref{eq:for5}) implies:
\be
E_D(T_\umu)\kformep=\sum_{\ulamb_{\umu}\in {\cal P}_h}
b_{\ulamb_{\umu}}E_D(\Delta_{\ulamb_{\umu}}(\TT))\kformep\label{eq:for6}
\ee
The left hand side is a $\ZZ$-linear combination of $\ikformep$ equal to:
\[
\sum_{\ulamb_{\umu}\in {\cal
P}_h}b_{\ulamb_{\umu}}\epsilon^{1+r_1(\ulamb_{\umu})}
\wedge\ldots\wedge \epsilon^{k+r_k(\ulamb_{\umu})}.
\]
By unicity, the coefficients $b_{\ulamb_{\umu}}$'s must be the same as those
occurring in the expansion of the l.h.s. of~(\ref{eq:for6}) in terms of the
$\ikformep$, i.e. they must be integers.
\qed

\begin{prop}\label{prop3.6}
Let $\phi\in\ZZ[T_1,\ldots,T_k]\subseteq \ZZ[\TT]$ such
that
\[
\phi(D_1,\ldots,D_k)\kformep=0,
\]
then $\phi=0$.
\end{prop}
\proof
Let $\phi\in \ZZ[T_1,\ldots,T_k]$ be a homogeneous polynomial of degree $h$
such that $\phi(D_1,\ldots, D_k)\kformep=0$. Since
$\{\Delta_\ulamb(\TT):|\ulamb|=h\}$, is a
$\ZZ$-basis of $\ZZ[\TT]_h$, then:
\[
\phi(T_1,\ldots, T_k)=\sum_{\ulamb\in{\cal P}_h}a_\ulamb\Delta_\ulamb(\TT).
\]
Hence:
\[
0=\phi(D_1,\ldots, D_k)\kformep=\sum_{\ulamb\in{\cal
P}_h}a_\ulamb\Delta_\ulamb(D)\kformep,
\]
implying that all $a_\ulamb=0$, because $\Delta_\ulamb(D)\kformep$ are linearly
independent.

\qed

Given that $\ker(E_{\kformep})$ is a homogeneous ideal, the grading of $\ZZ[D]$
induces the {\em codimension grading} of
\[
{\cal A}^*(\bigwedge^kM,D_t):={\ZZ[D]\over\ker(E_{\kformep})},
\]
the {\em intersection ring} of the pair $(\bigwedge^kM,D_t)$, namely:
\[
{\cal A}^*(\bigwedge^kM,D_t)={\cal A}^0(\bigwedge^kM,D_t)\oplus{\cal
A}^1(\bigwedge^kM,D_t)\oplus{\cal A}^2(\bigwedge^kM,D_t)\oplus\ldots
\]
The {\em codimension $i$} submodule ${\cal
A}^i(\bigwedge^kM,D_t)$ of ${\cal A}^*(\wkM,D_t)$ is, by definition, 
the module generated by the classes
modulo
$\ker(E_{\kformep})$ of the monomials in $\ZZ[D]$ of degree $i$, with respect
to the grading~(\ref{eq:forgrading}).
It turns out that $\bigwedge^kM$ is a  graded (by codimension) ${\cal
A}^*(\bigwedge^kM,D)$ module, i.e.:
\[
\bigwedge^kM=\ZZ\oplus (\bigwedge^kM)_1\oplus (\bigwedge^kM)_2\oplus\ldots
\]
Clearly, $(\bigwedge^kM)_h$ is the (isomorphic) image of ${\cal
A}^h(\bigwedge^kM,D)$ through the evaluation morphism: it is spanned by all the
$k$-vectors $\ikformep$ such that $\sum_{j=1}^k i_j=1+2+\ldots+k+h=k(k+1)/2+h$,
i.e. by all the $k$-vectors $\ikformep$ of {\em weight} $h$, where the
{\em weight} of a $k$-vector is, by definition, $\sum_{j=1}^k(i_j-j)$.

\subsubsection{The Intersection Ring}

\begin{lemma}\label{leminv}
Let $E_t$ be the inverse of the ${\cal S}$-derivation $D_t$ (see
Section~\ref{sect1.2}). Then for all
$k>0$ and all
$h>k$,
$E_h\alpha=0$ for all $\alpha\in \wkM$.
\end{lemma}
\proof
The proof works by induction on $k$. The property is true for $k=1$, because
$E_2=D_1^2-D_2=0$ in
$M$. Suppose it
holds for $k-1$. If $\alpha\in \wkM$,  by Theorem~\ref{giamth}, there is
$G_\alpha\in \ZZ[\TT]$, such that
  $G_\alpha(D)\kformep=\alpha$. One may thence check the property for
$\alpha=\kformep$, since
\[
E_h(\alpha)=E_h\circ
G_\alpha(D)\kformep=G_\alpha(D)E_h(\kformep).
\]
Let
$\beta=\epsilon^2\wedge\ldots\wedge\epsilon^k$. Then:
\[
E_h(\epsilon^1\wedge\beta)=\sum_{i=0}^{k-1}E_{h-i}(\epsilon^1)\wedge
E_i(\beta)+E_1(\epsilon^1)\wedge E_{h-1}\beta+\epsilon^1\wedge E_h(\beta)=0,
\]
and this proves the claim.
\qed

Let ${\cal A}^*(\wkM,D)$ be the intersection ring of $\wkM$.
\begin{prop}\label{propgir}
The ring ${\cal A}^*(\wkM,D)$ is the polynomial ring $\ZZ[D_1,\ldots,
D_k]$.
\end{prop}
\proof
By Proposition~\ref{prop3.6}, $D_1,\ldots, D_k$ are algebraically independent.
Now we shall show that $D_{k+i}$ ($i\geq 1$) is a $\ZZ$-polynomial
in $D_1,\ldots,D_k$.
In fact, for each
$j\geq 0$, relation $E_tD_t=D_tE_t=1$ implies:
\be
D_{k+i}-E_1D_{k+i-1}+\ldots+(-1)^{k+i-1}E_{k+i-1}D_1+
(-1)^{k+i}E_{k+i}=0.\label{eq:fordip}
\ee
Since $E_{p}$ is a polynomial expression in $D_1,\ldots, D_p$ only, it follows
that (for $i=1$) $D_{k+1}$ is a polynomial in $D_1,\ldots, D_k$, because
$E_{k+1}=0$ by Lemma~\ref{leminv}. Suppose now that the property holds for all
$1\leq j\leq i-1$. Then formula~(\ref{eq:fordip}) reads as:
\be
D_{k+i}=-E_1D_{k+i-1}+\ldots+(-1)^{k+1}E_{k}D_i,\label{eq:polrel}
\ee
because $E_{k+i}=0$ for all $i\geq 1$; it follows that $D_{k+i}$ is a 
polynomial
in
$D_1,
\ldots, D_k$, by induction. As a conclusion:
\[
{\cal A}^*(\wkM,D_t)={\ZZ[D_1,D_2,\ldots]\over
(D_{k+i}-D_{k+i}(D_1,\ldots,D_k))_{i\geq 1}}\cong \ZZ[D_1,\ldots, D_k].
\]
\qed

\begin{rmk}
 From now on, when working on $\wkM$, the notation $D_{h}$
for
$h>k$, will stand for the polynomial expression in
$D_1,\ldots, D_k$ implicitly defined by formula~{\rm (\ref{eq:polrel})}.
\end{rmk}

Let $D_nM$ be the $\ZZ$-submodule of $M$ generated by
$(\epsilon^{n+1}, \epsilon^{n+2},\ldots,)$ and let 
$\bigwedge^{k-1}M\wedge D_nM$
be the
$\ZZ$-submodule of $\wkM$ generated by all the $k$-vectors $\ikformep$ with
$i_k\geq n+1$.
\begin{prop}\label{propideale}
Let $I=(i_1,\ldots,i_k)$ be a symbol (Sect.~\ref{sect1.1}) and suppose that
\[
\Delta_{\ulamb(I)}(D)\kformep=\ikformep\in \bigwedge^{k-1}M\wedge D_nM.
\]
Then $\Delta_{\ulamb(I)}(D)$ belongs to the ideal of $\ZZ[D]$ generated by
$(D_{n-k+1},\ldots,D_n)$.
\end{prop}
\proof
If $k+r_k=i_k\geq n+1$, one has:
\[
\Delta_{\ulamb(I)}=\left|\matrix{D_{r_1}&D_{r_2+1}&
\ldots&D_{r_k+k-1}\cr
D_{r_1-1}&D_{r_2}&\ldots&D_{r_k+k-2}\cr
\vdots&\vdots&\ddots&\vdots\cr
D_{r_1-k}&D_{r_2-k+1}&\ldots&D_{r_k}}\right|,
\]
which proves, by expanding the determinant along last column, that
$\Delta_{\ulamb(I)}\in (D_{n-k+1}, D_{n-k+2},\ldots)$.

It suffices then to
prove that for $i\geq 1$, $D_{n+i}\in (D_{n-k+1}, D_{n-k+2},\ldots, D_n)$ .
Observe now that:
\[
D_{n+1}=\sum_{j=1}^{k}(-1)^{j-1}E_jD_{n+1-j}+\sum_{j=k+1}^{n+1}
(-1)^{j-1}E_iD_{n+1-j}=\sum_{j=1}^{k}(-1)^{j-1}E_jD_{n+1-j}
\]
because $E_{h}(D)=0$ if $h>k$. This proves that $D_{n+1-j}\in (D_{n},\ldots,
D_{n-k+1})$ and hence, by induction on $i$, that
\[
D_{n+i}\in (D_{n},\ldots,
D_{n-k+1}),\quad \forall h>n.
\]

\qed

\section{More Intersection Pairs}\label{Sect4}
\subsection{The intersection pair $\left(\bigwedge^k(M/D_nM), D_t\right)$}

Let $(M,D_t)$ be the canonical intersection pair (Sect.~\ref{Sect3.1}) and let
$M_n=M/D_nM$. It is a
  $\ZZ$-module freely generated by $\epsilon^i+D_nM$, which, by abuse, will be
indicated as $(\epsilon^1, \ldots, \epsilon^n)$. The map $D_t:M\lra M[[t]]$
induces a map $D_t:M_n\lra M_n[[t]]$ in a obvious way.

\begin{prop}\label{presip}
For each $k\geq 1$, $(\wkM_n,D_t)$ is an ${\cal IP}$. Its
intersection ring is:
\[
{\cal A}^*(\wkM_n,D)={\ZZ[D_1,\ldots,D_k]\over(D_{n-k+1},\ldots, D_n)}
\]
\end{prop}
\proof
In fact:
\[
\wkM_n=\bigwedge^k \left({M\over D_nM}\right)\cong{\wkM\over
\bigwedge^{k-1}M\wedge D_nM}
\]
Clearly there is an epimorphism $\ZZ[D_1,\ldots, D_k]\lra \wkM_n$ gotten
by composing the isomorphism $\ZZ[D_1,\ldots, D_k]\cong
\wkM$ with the canonical projection of $\wkM$ onto $\wkM_n$. We claim that:
\[
D_{n-k+i}(\ikformep)=0\,\,{\rm mod}\,\, \bigwedge^{k-1}M\wedge D_nM.
\]
This is because (notation as in Sect.~\ref{sect1.1})
\begin{eqnarray*}
&{}&D_{n-k+j}(\ikformep)=D_{n-k+i}\Delta_{\ulamb(I)}(D)\kformep=\\
&{=}&\Delta_{\ulamb(I)}(D)D_{n-k+i}\kformep=\quad (by\,\, Pieri's\,\,
formula)=\\ &=&\Delta_{\ulamb(I)}(D)(\epsilon^1\wedge\ldots
\epsilon^{k-1}\wedge \epsilon^{n+j})=0\,\,{\rm mod}\,\,
\bigwedge^{k-1}M\wedge D_nM ,
\end{eqnarray*}
proving that
${\ZZ[D_1,\ldots,D_k]\over(D_{n-k+1},\ldots, D_n)}$ surjects onto $\wkM_n$.
There are no relations in degree less than $n-k$, because any such would
also be a relation in ${\cal A}^*(\wkM,D)$, contradicting
Proposition~\ref{prop3.6}.

Suppose there is a relation of degree  $n-k+1+j$ ($j\geq 0$), i.e.
\be
\phi(D_1,\ldots, D_k)\kformep\in \bigwedge^{k-1}M\wedge D_nM.\label{eq:correl}
\ee
By Lemma~\ref{deltabasis}, such a $\phi$ can be written in a 
unique way as:
\[
\phi(D_1,\ldots,D_k)=\sum 
a_{r_1,\ldots,r_k}\Delta_{r_k,\ldots,r_1}(D),
\]
for some integers 
$a_{r_1,\ldots,r_k}$.
Equation~(\ref{eq:correl}) is thence equivalent 
to:
\[
\sum a_{r_1\ldots
r_k}\Delta_{r_k,\ldots,r_1}(D)\kformep=a_{r_1\ldots
r_k}\rikformep\in \bigwedge^{k-1}M\wedge D_nM.
\]

By Proposition~\ref{propideale}, each  $\Delta_{(r_k\ldots 
r_1)}(D)\in (D_{n-k+1},
D_{n-k+2},...,D_n)$ and hence $\phi(D_1,\ldots,D_k)\in (D_{n-k+1},
D_{n-k+2},...,D_n)$, too.
\qed
\begin{corol}\label{corisoclas}
There is a canonical ring isomorphism  
\[
\Phi: {\cal
A}^*(\wkM_n,D)\lra A^*(G),
\] 
where $A^*(G)$ is the Chow ring of 

$G:=G_k(V)$, $\dim(V)=n$.
\end{corol}
\proof
Let ${\cal A}^*$ be the ring:
\[
{\cal A}^*:={\ZZ[D_1,\ldots, D_{n-k}]\over(Y_{k+1}(D),\ldots,Y_n(D))}
\]
where $Y_i(D)$'s are implicitly defined by:
\[
{1\over 1+D_1t+\ldots+D_{n-k}t^{n-k}}=\sum_{i\geq 0}(-1)^iY_i(D)t^i,
\]
We then show that  ${\cal A}^*\cong{\cal A}^*(\wkM_n,D_t)$.
Suppose first that $k<n-k$. Then by the very definition of $E_i$ one has that
$Y_j(D)=E_j(D)$ for all $1\leq j\leq n-k$. Since $E_{k+h}(D)=0$ for all $h\geq
1$, one has that $Y_{k+1}(D)=\ldots=Y_{n-k}(D)=0$ is equivalent to say that
$D_{k+1},\ldots,D_{n-k}$ can be expressed as a polynomial expression of
$D_1,\ldots, D_k$ only. Hence the ring ${\cal A}^*$ is indeed isomorphic to
\[
\ZZ[D_1,\ldots,D_k]\over(Y_{n-k+1}(D),\ldots,Y_n(D))
\]
Now, by looking at the definition of the $E_i$'s and of the $Y_i$'s 
one sees that
\be
E_{n-k+j}(D)=(-1)^{n-k+j-1}D_{n-k+j}+Y_{n-k+j}(D)\quad \forall j\geq
1\label{eq:forcoroll}
\ee
But $E_{n-k+j}(D)=0$, for all $j$ (because $n-k>k$), hence $Y_{n-k+j}(D)=0\iff
D_{n-k+j}=0$. Therefore ${\cal A}^*$ is indeed isomorphic to
\[
{\ZZ[D_1,\ldots,D_k]\over(D_{n-k+1},\ldots,D_n)}={\cal A}^*(\wkM_n,D_t).
\]
The case $n-k\leq k$ can be verified by using the obvious isomorphism between
$\bigwedge^{n-k}M_n$ and $\wkM_n$. However we present a direct check. One has:
\[
{\ZZ[D_1,\ldots,D_{n-k},D_{n-k+1},\ldots,D_k]\over
(D_{n-k+1},\ldots,D_k,D_{k+1},\ldots,D_n)}\cong
{\ZZ[D_1,\ldots,D_{n-k}]
\over (D_{k+1},\ldots,D_n)},
\]
but, again, the polynomial expression of  $D_{k+1},\ldots, D_n$ as a function
of $D_1, D_{n-k}$ is precisely $Y_{k+1}(D),\ldots,Y_n(D)$, so that:
\[
{\cal A}^*(\wkM_n,D_t)={\ZZ[D_1,\ldots,D_{n-k}]\over(Y_{k+1}(D),\ldots,
Y_n(D))}.
\]
It is then clear that the map $\Phi: {\cal A}^*(\wkM,D_t)\lra A^*(G)$ 
defined by
$D_i\mapsto\sigma_i$ is an isomorphism, because of the
presentation~(\ref{eq:prescrg}).

\qed

\subsection{(Small) Quantum Intersection Pair}
Let $(M,D_t)$ be the canonical intersection pair and let
$M[D_n]:=M\otimes_{\bf Z}\ZZ[D_n]$, $n\geq 1$. The module $M[D_n]$ is nothing
but
$M$ itself thought of as a $\ZZ[D_n]$-module. It is a finitely generated free
$\ZZ[D_n]$-module of rank $n$, generated by $(\epsilon^1,\ldots, \epsilon^n)$.
Denote by the same symbol the
$D_n$-linear extension of $D_t$
\begin{prop}
The pair $(M[D_n],D_t)$ is a $\ZZ[D_n]$-intersection pair. Its intersection
ring is:
\be
{\cal A}^*(M[D_n],D_t)={\ZZ[D_n][D_1]\over (D_1^n-D_n)}.
\ee
\end{prop}
\proof
  The natural evaluation map $\ZZ[D_n][\TT]\lra M[D_n]$ is clearly surjective,
because $M\cong\ZZ[D_1]$. It is only left to compute the intersection ring. One
has that $D_i=D_1^i$. Hence it is generated by $D_1$ only. Furthermore one has
the relation $D_1^n=D_n$, whence the claim
\qed

By renaming $D_n$ with $q$, one may write
\[
{\ZZ[D_n,D_1]\over (D_1^n-D_n)}\cong {\ZZ[q,D_1]\over (D_1^n-D_n, D_n-q)}\cong
{\ZZ[q,D_1]\over (D_1^n-q)}.
\]
Last expression is exactly that of the {\em small quantum cohomology ring of}
$\PP^{n-1}$ (see e.g.~\cite{FuPan}), once $D_1$ is interpreted as the 
class of the hyperplane section. For this reason the pair
$(M[q],D_t)$ shall be referred to as to the {\em quantum canonical
intersection pair}. Let us now extend $D_t$ to an $HS$-derivation on $\wkM[q]$,
thought of as a
$\ZZ[q]$-module. Still one has:
\begin{prop}
For each $1\leq k\leq n$, $(\bigwedge^kM[D_n], D_t)$ is a
$\ZZ[D_n]$ cohomological pair generated by $\kformep$.
\end{prop}
  \proof It is sufficient to show that for each $\ikformep$ there exists
$P(\TT)\in\ZZ[D_n][\TT]$ such that $E_D(P(\TT))\kformep=\ikformep$. This is
a consequence of Giambelli's formula, because
$\Delta_\ulamb(D)\in\ZZ[D]=\ZZ[D_n][D]$. In other words:
\[
\ikformep=\Delta_\ulamb(D)\kformep.
\]

\begin{prop}
The intersection ring of the $\ZZ[D_n]-{\cal IP}$ $\left(M[D_n]\right),D_t)$ is
given by
\[
{\cal
A}^*(\bigwedge^k\left(M[D_n]\right),D_t)={\ZZ[D_n][D_1,\ldots,D_k]
\over(D_{n-k+1},\ldots,D_{n-1})}
\]
\end{prop}
\proof
Clearly $(D_{n-k+1},\ldots, D_{n-1})$ is contained in the
ideal ${\cal R}$ of relations. In fact one has, for $1\leq i\leq k-1$
\begin{eqnarray*}
D_{n-k+i}\kformep&=&
\epsilon^1\wedge\ldots\wedge\epsilon^{k-1}\wedge\epsilon^{n+i}=\\
\epsilon^1\wedge\ldots\wedge\epsilon^{k-1}\wedge
D_n\epsilon^{i}&=&D_n(\epsilon^1\wedge\ldots\wedge\epsilon^
i\wedge\ldots\wedge\epsilon^{k-1}
\wedge\epsilon^i)=0,
\end{eqnarray*}
where the last vanishing is due to skew-symmetry.
We contend that indeed ${\cal R}=(D_{n-k+1},\ldots,D_n)$. We first
observe that over the integers there is no relation of degree $\leq 
n-k$, because otherwise it would also be a relation in ${\cal 
A}^*(\wkM,D_t)$.
  Suppose that $0\neq \psi\in
\ZZ[D_n][T_1,\ldots,T_k]$ is such that $\psi(D)\kformep=0$ and let
$h=\deg(\psi)$. Denote by $\tilde{\psi}(D)$ the same $\psi$ thought of as a
polynomial with integral coefficients.
If
\[
\tilde{\psi}(D)\kformep=0,
\]
then necessarily $\tilde{\psi}(D)\kformep\in\bigwedge^{k-1}M\wedge D_nM$.
It follows (Prop.~\ref{propideale}) that $\tilde{\psi}(D)\in (D_{n-k+1},\ldots,
D_n)$, that is
$\tilde{\psi}(D)\kformep=\chi(D)D_n\kformep$, where $\chi$ is a 
$\ZZ$-polynomial
of degree
$h-n$, such that
$\psi(D)\kformep=0$. Hence, by repeating the argument, one would conclude the
existence of a non zero relation over the integers of degree $\leq n-k$, and
this is not possible, by Prop.~\ref{prop3.6}.
\qed

By introducing a new indeterminate $q$, the intersection ring can be expressed
as:
\[
{\cal
A}^*(\bigwedge^k M[D_n],D_t)={\ZZ[D_n][D_1,\ldots,D_k]
\over(D_{n-k+1},\ldots,D_{n-1})}={\ZZ[q][D_1,\ldots,D_k]
\over(D_{n-k+1},\ldots,D_{n-1}, D_n-q)}.
\]

\begin{corol}\label{corwstder}
There is a canonical ring isomorphism
\[
\Phi_q:  {\cal A}^*(\wkM_n[q],D)\lra QA^*(G),
\]
  where $QA^*(G)$ is the quantum Chow ring of  $G:=G_k(V)$, $\dim(V)=n$.
\end{corol}
\proof
As in Corollary~\ref{corisoclas} one starts by observing the
isomorphism:
\be
{\ZZ[q,D_1,\ldots,D_{n-k}]\over(Y_{k+1}(D),\ldots,Y_n(D)-(-1)^{n-k-1}q)}\cong
{\ZZ[q,D_1,\ldots,D_{k}]\over(D_{n-k+1},\ldots,D_n-q)}.\label{eq:isomanc}
\ee
By formula~(\ref{eq:forcoroll}) and the fact that $E_{k+j}(D)=0$ for all $j\geq
1$, one has that $(-1)^{h}D_h(D_1,\ldots,D_k)=Y_h(D)$ for all $h\geq k+1$. In
particular
$(-1)^{n}D_n=Y_n(D)=(-1)^{n-k+1}q$. It follows that the relation
$Y_n(D)-(-1)^{n-k-1}q$ is equivalent to $D_n-(-1)^{k-1}q$. Therefore the
isomorphism~(\ref{eq:isomanc}) is gotten by $q\mapsto (-1)^{k-1}q$.
The conclusion is that the morphism
\[
{\cal
A}^*(\wkM_n[q],D)\lra \Phi_q:QA^*(G),
\]
sending $D_i\mapsto \sigma_i$ and $q\mapsto (-1)^{k-1}q$ is
an isomorphism, because of the presentation~(\ref{eq:wstring}), as claimed.

\qed
\begin{prop}\label{Pieri} For each $1\leq h\leq k$, {\em Pieri's formula for
quantum
${\cal S}$-derivatives} holds:
\be
D_h(\ikformep)=\sum_{J}\epsilon^{j_1}\wedge\ldots\wedge
\epsilon^{j_k}+(-1)^{k-1}q
\sum_{J'}\epsilon^{j'_1}\wedge\ldots\wedge \epsilon^{j'_k},\label{eq:qpier}
\ee
where the first sum is over all $J=(j_1,\ldots,j_k)$ such that
\[
i_1\leq j_1<
i_2\leq j_2<\ldots< i_k\leq j_k\leq n, \quad
|J|=|I|+h,
\]
and the second sum
is over all $J'=(j'_1,\ldots,j'_k)$ such that $|J'|=|I|+h-n$ and
\[
1\leq j'_1<i_2\leq
j_2'\leq i_3<\ldots<i_{k-1}\leq j_{k-1}'<i_k.
\]
\end{prop}
\proof

First of all observe that if $q^m\epsilon^i$ is written as $\epsilon^{mn+i}$,
then, for each
\[
1\leq i_1<\ldots<i_k\leq n
\]
one has, by Pieri's for ${\cal
S}$-derivations:
\[
D_h(\ikformep)=\sum\epsilon^{i_1+h_1}\wedge\ldots\wedge\epsilon^{i_k+h_k}
\]
where the sum is over all non-negative $(h_1,\ldots,h_k)$ such that:
\[
i_1\leq i_1+h_1<i_2\leq i_2+h_2<\ldots <i_{k-1}\leq
i_{k-1}+h_{k-1}< i_k\leq i_k+h_k.
\]
Last sum may be split as:
\[
\sum_{i_k+h_k\leq n }\ihikformep+\sum_{i_k+h_k>n}\ihikformep\quad (\sum h_p=h)
\]
where in both cases $i_1+h_1<\ldots<i_k+h_k$ and $\sum_{p=1}^kh_p=h$.
The second summand may hence be written as:
\begin{eqnarray}
&{ }&\sum_{i_k+h_k>n}\ihikformep=\sum_{i_k+h_k>n}\epsilon^{i_1+h_1}\wedge\ldots\wedge 
\epsilon^{i_{k-1}+h_{k-1}}\wedge 
q\epsilon^{i_k+h_k-n}\nonumber\\&{}&=(-1)^{k-1}q\sum_{i_k+h_k>n}
\epsilon^{i_k+h_k-n}\wedge\ldots\wedge\epsilon^{i_j+h_j}
\wedge\ldots\wedge
\epsilon^{i_{k-1}+h_{k-1}}\label{eq:uglex}
\end{eqnarray}
Now, if $h_k$ is such that $i_k-n+h_k>i_j+h_j$ ($2\leq j\leq k-2$), there are
$h'_j$
$h'_k$ such that
$h'_j+h'_k=h_j+h_k$, $i_k-n+h'_k=i_j+h_j$ and $i_j+h'_j=i_k-n+h_k$. It
suffices to choose $h'_j= i_k-n+h_k-i_j$ and $h'_k=i_j+h_j+n-i_k$ (obviously
$h'_j\geq 0$ and $h'_k\geq 0$ since $i_k\leq n$). It follows that
in  sum~(\ref{eq:uglex}) each term corresponding to $h_k$ such that
$i_k-n+h_k>i_j+h_j$ cancel out against a term of the form:
\[
(-1)^{k-1}q\epsilon^{i_k+h'_k-n}\wedge\epsilon^{i_2+h_2}
\wedge\ldots\wedge \epsilon^{i_j+h'_j}\wedge\ldots\wedge
\epsilon^{i_{k-1}+h_{k-1}}.
\]
Therefore, sum~(\ref{eq:uglex}) is the same as:
\[
(-1)^{k-1}q\sum_{i_k+h_k>n}
\epsilon^{i_k+h_k-n}\wedge\ldots\wedge\epsilon^{i_j+h_j}
\wedge\ldots\wedge\epsilon^{i_{k-1}+h_{k-1}},
\]
such that $i_k+h_k-n<i_1+h_1<i_2+h_2<\ldots<i_{k-1}+h_{k-1}$. Putting
$j'_1=i_k+h_k-n$ and $j'_p=i_p+h_p$, one has exactly formula~(\ref{eq:qpier}).

\qed

\section{Geometry}\label{Sect5}
\subsection{ }

The purpose of this concluding section is to show how (quantum) Schubert
calculus for Grassmannians  fits into the algebraic framework studied in the
previous sections.
Let $G:=G_k(V)$ be the grassmannian of $k$-planes in a $n$-dimensional
vectorspace $V$ over an algebraically closed field $\kk$. A $k$-plane
$[\Lambda]$ is determined, modulo the natural action of the group $Gl_k(\kk)$,
by a ``column"
\[
\Lambda=\pmatrix{v_1\cr\vdots\cr v_k}\in V^k
\]
of linearly independent vectors $v_i\in V$. If $\alpha\in V^\vee$, set:
\[
\alpha(\Lambda)=\pmatrix{\alpha(v_1)\cr\vdots\cr\alpha(v_k)}\in \kk^k.
\]
Therefore, if $(e_i)$ is a basis of $V$ and $(\epsilon^j)$ its dual:
\[
(\epsilon^1,\ldots,\epsilon^n)(\Lambda)=\left(\epsilon^1(\Lambda),\epsilon^2(\Lambda),\ldots,
\epsilon^n(\Lambda)\right)=
\pmatrix{\epsilon^1(v_1)&\epsilon^2(v_1)&\ldots&\epsilon^n(v_1)\cr
\vdots&\vdots&\ddots&\vdots\cr
\epsilon^1(v_k)&\epsilon^2(v_k)&\ldots&\epsilon^n(v_k)}.
\]
Similarly:
\[
\ikformep(\Lambda)=\epsilon^{i_1}(\Lambda)\wedge\ldots\wedge\epsilon^{i_k}(\Lambda)=
\left|\matrix{\epsilon^{i_1}(v_1)&\ldots&\epsilon^{i_k}(v_1)\cr
\vdots&\ddots&\vdots\cr
\epsilon^{i_1}(v_1)&\ldots&\epsilon^{i_k}(v_k)}\right|.
\]
For each symbol  $I=(i_1,\ldots,i_k)$ (see Sect.~\ref{sect1.1}) with $1\leq
i_1<i_2<\ldots<i_k\leq n$, there is an open set of the grassmannian
\[
G_I:=\{[\Lambda]\in G, \ikformep(\Lambda)\neq 0\}.
\]
The definition does not depend on the particular representative of $[\Lambda]$.
If $[\Lambda]\in G_I$, then $[\Lambda]=[\Lambda_I]$ where:
\[
\Lambda_I=(\epsilon^{i_1}(\Lambda),\ldots,\epsilon^{i_k}(\Lambda))^{-1}\Lambda.
\]
If $[\Lambda]\in G_I\cap G_J$, then
$
\Lambda_I=m_{IJ}\Lambda_J
$, where
\[
m_{IJ}=(\epsilon^{i_1}
(\Lambda),\ldots,\epsilon^{i_k}(\Lambda))^{-1}\cdot(\epsilon^{j_1}
(\Lambda),\ldots,\epsilon^{j_k}(\Lambda))\in Gl_k(\kk).
\]
  Clearly
$m_{II}=id_{\kk^k}$, and $m_{IJ}\cdot m_{JK}=m_{IK}$.
The rank $k$ vector bundle gotten by gluing all the products
$G_I\times
\kk^k$ via the $m_{IJ}$, is the tautological bundle ${\cal T}$. It 
sits into the
universal exact sequence:
\[
0\lra {\cal T}\lra O_{G}^{\oplus n}\lra {\cal Q}\lra 0,
\]
where ${\cal Q}$ is the {\em universal quotient bundle}. The bundle ${\cal T}$
can be also described as the set of pairs $([\Lambda],v)\in G\times \CC^n$ such
that $v\in [\Lambda]$.
Any linear form $\alpha\in V^\vee$ can be seen then as a map of
vector bundles:
\[
\alpha:{\cal T}\lra O_G.
\]
Similarly, any $k$-tuple $\alpha_1,\ldots,\alpha_k\in{\cal T}^\vee$  induces
a map
\[
\alpha_1\wedge\ldots\wedge\alpha_k:\wedge^k{\cal T}\lra O_G,
\]
i.e. a
global section of the bundle $O_G(1):=\bigwedge^k{\cal T}^\vee$. Clearly,
$
H^0(G,O_G(1))\cong \bigwedge^kV^\vee,
$
generated by the global sections
$\left\{\ikformep\right\}_{1\leq i_1<i_2<\ldots<i_k\leq n}$ .
The map
\[
\left\{\matrix{p&:&G&\lra&  \PP\left(\bigwedge^kV^\vee\right)\cong
\PP^{{n\choose k}-1}\,\,\,\,\,\,\,\,\,\,\,\,\,\,\,&{}\cr
{ }&{ }&[\Lambda]&\longmapsto&\,\,\,\,\,\,\,[\ikformep(\Lambda)]_{1\leq
i_1<\ldots<i_k\leq n}&{ }
}\right.,
\]
is the well known {\em
Pl\"ucker embedding} and the $\ikformep$ are the Plucker coordinates. We shall
see in a moment how Schubert varieties can be written in terms of
Pl\"ucker coordinates.
\subsection{ }
Let
\be
{\cal E}:V=E_0\supset E_1\supset\ldots\supset E_n=(0)\label{eq:flag}
\ee
be a complete flag of $V$, such that ${\rm codim}\,{E_i}=i$.
Let $(e_1,\ldots, e_n)$ be an ${\cal E}$-adpated basis, that is:
$E_i=span\{e_{i+1},\ldots, e_n\}$. In terms of the dual basis:
$E_i^\circ=span\{\epsilon^1,\ldots,\epsilon^i\}$, where $E_i^\circ\subseteq
V^\vee$ is the {\em annihilator} of $E_i$.

\begin{prop}
The following formula holds:
\[
\dim([\Lambda]\cap E_i)= k-{\rm rk}
(\epsilon^{1}(\Lambda),\ldots,\epsilon^i(\Lambda)).
\]
\end{prop}
\proof
Let $(e_i)$ be an ${\cal E}$-adapted basis of $V$ and let $(\epsilon^j)$ be the
dual.
Clearly $[\Lambda]\cap E_i$ is the locus of vectors $\vv\in [\Lambda]$ such
that
$\epsilon^1(\vv)=\ldots=\epsilon^i(\vv)=0$. The number of
independent equations is the codimension of $[\Lambda]\cap E_i$, which is hence
the rank of the matrix $(\epsilon^1,\ldots,\epsilon^i)(\Lambda)$.
\qed

The chain of inequalities~(\ref{eq:flagdim}) may be thence written as:
\[
0\leq {\rm rk}(\epsilon^1(\Lambda))\leq {\rm
rk}(\epsilon^1(\Lambda),\epsilon^2(\Lambda))\leq\ldots\leq
{\rm rk}(\epsilon^1(\Lambda),\epsilon^2(\Lambda),\ldots,\epsilon^n(\Lambda))=k,
\]
indicating that ${\rm
rk}(\epsilon^{1}(\Lambda),\ldots,\epsilon^{h}(\Lambda))$ is equal to the
numbers of jumps $\leq h$. This simple observation gives us a useful
characterization of a Schubert symbol, namely that
$I([\Lambda])=(i_1,\ldots,i_k)$ if and only if:
\begin{eqnarray}
i_1&=&\min\{h\geq
1:\epsilon^h(\Lambda)\neq 0\}\quad {\rm and}\nonumber\\
i_j&=&\min\left\{h>i_{j-1}:\epsilon^{i_1}(\Lambda)\wedge\ldots\wedge
\epsilon^{i_{j-1}}(\Lambda)\wedge\epsilon^h(\Lambda)\neq 0\right\} \quad
(1<j\leq k).\nonumber\\  &{ }&
\end{eqnarray}
The latter can be refined into the following:
\begin{prop}\label{propeqschub}
\be
I([\Lambda])=(i_1,\ldots,i_k)\iff\left\{
\matrix{&\jkformep(\Lambda)=0,&\quad\forall (j_1,\ldots,j_k)\prec
I([\Lambda])\cr &\ikformep(\Lambda)\neq 0&{ }}\right.\label{eq:defsc}
\ee
\end{prop}
\proof
If $I([\Lambda])=(i_1,\ldots, i_k)$, one already knows that
$\ikformep(\Lambda)\neq 0$. Suppose now that $(j_1,\ldots, j_k)\prec
(i_1,\ldots,i_k)$ such that $\jkformep\neq 0$. Let $h$ the least integer $\geq
1$ such that
$j_h<i_h$. Then
$\epsilon^{i_1}\wedge\ldots\wedge\epsilon^{i_{h-1}}\wedge \epsilon^{j_h}\neq
0$, contradicting the minimality of $i_h$. This proves one implication.

On the other hand, suppose that equations~(\ref{eq:defsc}) are fulfilled.
First of all let us show that $i_1$ is  the first Schubert jump. Suppose that
$\epsilon^{j_1}(\Lambda)\neq 0$ for $j_1<i_1$. Since ${\rm
rk}(\epsilon^1,\ldots,\epsilon^{i_2})\geq 2$, there would exist $j_2\leq i_2$
such that
$\epsilon^{j_1}\wedge\epsilon^{j_2}(\Lambda)\neq 0$. Suppose to have 
constructed
$(j_2,\ldots, j_{l-1})\preceq (i_2,\ldots, i_{l-1})$ such that
\[
\epsilon^{j_1}\wedge\ldots\wedge\epsilon^{j_{l-1}}(\Lambda)\neq 0.
\]
  Then, since for all $h\geq 0$, ${\rm
rk}(\epsilon^1,\ldots,\epsilon^{i_l})\geq l$,  there exists $j_l\leq i_l$ such
that
\[
\epsilon^{j_1}\wedge\ldots\wedge\epsilon^{j_{l-1}}\wedge\epsilon^{j_l}\neq 0,
\]
  proving the existence of $(j_1,\ldots,j_k)\prec (i_1,\ldots, i_k)$
such that $\jkformep([\Lambda])\neq 0$, against the hypothesis. Therefore $i_1$
is the first Schubert jump. Now, suppose that
$i_1,\ldots, i_{l-1}$ are the first $l-1$ Schubert jumps; suppose there is
$i_{l-1}<j_l<i_l$ such that
$(\epsilon^{i_1}\wedge\ldots\wedge\epsilon^{i_{l-1}}\wedge\epsilon^{j_l})([\Lambda])
\neq 0$. Then, since ${\rm
rk}(\epsilon^1,\ldots,\epsilon^{i_{l+h}})\geq l+h$ ($h=1,\ldots, k-l$), there
exists $j_{l+h}\leq i_{l+h}$ such that
\[
\epsilon^{i_1}\wedge\ldots\wedge\epsilon^{i_{l-1}}\wedge\epsilon^{j_l}\wedge
\epsilon^{j_{l+1}}\wedge\ldots\wedge\epsilon^{j_k}([\Lambda])\neq 0,
\]
still contradicting the hypothesis.

\qed

We proved that equations~(\ref{eq:defsc}) are rank conditions of any matrix
\[
\Lambda=\pmatrix{\vv_1\cr\vdots\cr \vv_k}
\]
representing the $k$-plane $[\Lambda]$. In terms of an ${\cal
E}$-adapted basis $(e_i)$,  such a $[\Lambda]$ is represented by the matrix:
\be
\matrix{\Lambda_I\cr{ }\cr{
}}\matrix{=\cr{ }\cr{
}}\matrix{\pmatrix{0&\ldots&0&1&*&\ldots&*&0&*&\ldots&*&0&*&\ldots&*\cr
          0&\ldots&0&0&0&\ldots&0&1&*&\ldots&*&0&*&\ldots&*\cr
\vdots&\ddots&\vdots&\vdots&\vdots&\ddots&\vdots&\vdots&\vdots&\ddots
&\vdots&\vdots&\vdots&\ddots&\vdots\cr
0&\ldots&0&0&0&\ldots&0&0&0&\ldots&0&1&*&\ldots&*}\cr
\matrix{{ }&{ }&{ }&{ }&{ }&{ }&{ }&{ }&{ }&{ }&{ }&{}&{ }&{ }&{ }\cr
{ }&{ i_1}&{ }&&{ }&{ }&{ }&{i_2 }&{ }&{ }&{ }&{}&{ }&{ i_k}&{ }} }
\ee
showing that the set of all $k$-planes having $(i_1,\ldots,i_k)$ as Schubert
symbol, is an affine cell of codimension $(i_1-1)+\ldots +(i_k-k)$.

The Schubert cycle
$[\overline{W_{\ulamb(I)}({\cal E})}]$  will be denoted,
from now on, by the symbol:
$
[\ikformep]
$
meaning with such a notation the class in $A_*(G)$ of the closure of the
${\cal E}$-Schubert cell defined by equations~(\ref{eq:defsc}).
Let $M^\vee$ be the integral lattice of $V^\vee$ generated by the $\epsilon^i$.
Then there is a $\ZZ$-linear map:
\[
[\,\,\,]:\bigwedge M^\vee\lra A_*(G)
\]
defined by
\be
\sum_{1\leq i_1<\ldots < i_k\leq n}a_{i_1\ldots i_k}\ikformep\longmapsto
\sum_{1\leq i_1<\ldots < i_k\leq n}a_{i_1\ldots
i_k}[\ikformep].\label{eq:modiso}
\ee
This map would be an isomorphism were $A_*(G)$ torsion free, namely
were $[\ikformep]$ a $\ZZ$-basis for $A_*(G)$. This is well known (\cite{Fu1},
p.~27 and the references therein). However we can provide, within our
formalism, a new proof of this fact.
\begin{prop}\label{propfreelygen}
The classes $[\ikformep]$ freely generate $A_*(G)$.
\end{prop}
\proof

The map~(\ref{eq:modiso}) can be extended by ${\cal A}^*(\wkM,D_t)$-linearity
by defining
\[
P(D)[\ikformep]=[P(D)\ikformep].
\]
  If $\deg P(D)=j$, $P(D)A_i(X)\subseteq
A_{i-j}(X)$. This shows that any  $P(D)\in {\cal A}^j(\wkdM,D)$  corresponds to
a class in
$A^j(G)$.
Observe now that if $[\ikformep]\in A_i(G)$ then $D_1[\ikformep]\in A_{i-1}(G)$
and it is a $\ZZ$-linear combination of $[\jkformep]\in A_{i-1}(G)$ with
positive coefficients. This follows from the equality
$D_1[\ikformep]=[D_1(\ikformep)]$ and by Pieri's formula for ${\cal
S}$-derivatives (Prop.~\ref{Pieri}) applied to this special case, and the
morphism~(\ref{eq:modiso}).

Suppose now $[\ikformep]\in A_*(G)$ is such that $a[\ikformep]=0$
in $A_*(G)$, for some $0\neq a\in\ZZ$. Then:
\[
0=D_1^{k(n-k)-\sum
i_p}a[\ikformep]=ab[\epsilon^{n-k+1}\wedge\ldots\wedge\epsilon^n]
\]
for some positive integer $b$. In fact
\[
D_1^{k(n-k)-\sum
i_p}A_{|\ulamb|}(G)=A_{k(n-k)}(G)\cong
\ZZ\cdot[\epsilon^{n-k+1}\wedge\ldots\wedge\epsilon^n],
\]
where $[\epsilon^{n-k+1}\wedge\ldots\wedge\epsilon^n]$ is the class of the
point. Therefore $ab[\epsilon^{n-k+1}\wedge\ldots\wedge\epsilon^n]=0$ if and
only if $ab=0$, i.e. only if $a=0$, contradicting the assumption. It follows
that
$A_*(G)$ is a torsion free
$\ZZ$-module, hence free~(\cite{Lang}, p.~147).
\qed

\begin{prop}
Giambelli's formula holds:
\[
[\rikformep]=\Delta_\ulamb(D)[\kformep],
\]
where $\ulamb=(r_k\ldots r_1)$.
\end{prop}
\proof
In fact:
\[
[\rikformep]=[\Delta_\ulamb(D)\kformep]=\Delta_\ulamb(D)[\kformep]
\]
\qed

One concludes that $\Delta_\ulamb(D)\in A^*(G)$ is the Poincar\'e dual of
$[\rikformep]$. As remarked, the ring ${\cal A}^*(\wkM,D_t)$ is 
isomorphic, as a
$\ZZ$-module, to $A^*(G)$. The isomorphism $\Phi_q$ (Cf.
Corollary~\ref{corwstder}) sends
$\Delta_\ulamb(D)$ onto $\sigma_\ulamb$ and gives to $A_*(G)$ another structure
of ${\cal A}^*(\bigwedge^k\dM,D_t)$-module, namely:
\[
D_h*[\ikformep]=\sigma_h\cap[\ikformep].
\]

We contend that:
\[
D_h[\ikformep]=\sigma_h\cap[\ikformep].
\]
\begin{thm}\label{thmiso}
The two ${\cal
A}^*(\bigwedge^k\dM,D)$-module structures of $A_*(G)$ coincide,
i.e.
\[
D_h[\ikformep]=\sigma_h\cap[\ikformep].
\]
\end{thm}
\proof\,
By {\em Pieri's formula for ${\cal S}$-derivatives} (Prop.~\ref{Pieri}), one
has:
\[
D_h[\ikformep]=\sum_{\matrix{i_1\leq j_1<\ldots <
i_k\leq j_k\cr |J|=|I|+h}}[\epsilon^{j_1}\wedge\ldots\wedge
\epsilon^{j_k}],
\]
where $J=(j_1,\ldots, j_k)$, $I=(i_1,\ldots,i_k)$, $|J|=\sum_p j_p$ and
$|I|=\sum_p i_p$. Hence, in terms of
$\sigma$'s, one has:
\begin{eqnarray}
D_h[\ikformep]&=&\sum_{i_1\leq j_1<i_2\leq j_2<\ldots<i_k\leq j_k}
\sigma_{j_k-k,\ldots,j_1-1}\cap
[\kformep]\nonumber \\&=&\sum_{\scriptsize{\matrix{\lambda_1\geq
\mu_1\geq\ldots\geq\lambda_k\geq\mu_k\cr
|\umu|=|\ulamb|+h}}}{\sigma_\umu}\cap[\kformep],\label{eq:pir3}
\end{eqnarray}
having set $\mu_i=i_{k-i+1}-(k-i+1)$ and
$\lambda_i=i_{k-i+1}-(k-i+1)$. Furthermore, by Pieri's
formula~(\ref{eq:piericlas}):
\begin{eqnarray}
\sum_{\scriptsize{\matrix{\lambda_1\geq
\mu_1\geq\ldots\geq\lambda_k\geq\mu_k\cr
|\umu|=|\ulamb|+h}}}{\sigma_\umu}\cap[\kformep]&=&
(\sigma_h\cup\sigma_\ulamb)\cap
[\kformep]=\nonumber\\{ }\nonumber\\=\sigma_h\cap
(\sigma_\ulamb\cap[\kformep])&=& \sigma_h\cap
[\ikformep].\label{eq:pir4}
\end{eqnarray}
Therefore the first member of~(\ref{eq:pir3}) is equal to the second
member of~(\ref{eq:pir4}), as claimed.

\qed

We may hence conclude that Schubert Calculus for grassmannians is entirely
encoded by the algebra of the intersection pair
$(\bigwedge^k\dM,D_t)=\bigwedge^k(\dM,D_t)$, where $(\dM,D_t)$ is the canonical
${\cal IP}$.

\subsection{Quantum Schubert Calculus}

By Corollary~\ref{corwstder} there is  an isomorphism
$\Phi_q: {\cal A}^*(M[q],D)\lra QA^*(G)$, defined by $D_i\mapsto\sigma_i$
and
$q\mapsto (-1)^{k-1}q$. As for ``classical" Schubert calculus, this isomorphism
extends to an isomorphism of modules:
\[
\sigma_\ulamb\cap_q[\ikformep]=\Phi_q(\Delta_\ulamb(D)[\ikformep]).
\]
This can be achieved by observing that by Bertram-Pieri's formulas:
\[
\sigma_h\cap_q[\ikformep]=(\sigma_h\cup_q\sigma_\ulamb)\cap_q
[\kformep]=\sigma_h\cup\sigma_\ulamb +q\sum \sigma_\mu,
\]
where the sum is extended to all
\[
\lambda_1-1\geq\mu_1\geq \lambda_2-1\geq\ldots\geq \lambda_k-1\geq \mu_k\geq 0
\]
such that $|\nu|=|\mu|+h-n$.
On the other hand, by quantum Pieri's formula for Schubert derivations one has:
\[
D_h[\ikformep]= D_h[\ikformep] {\rm mod}\,\bigwedge^{k-1}M\wedge D_nM
+(-1)^{k-1}\sum\jkformep
\]
where the sum is extended to all
\[
1\leq j_1<i_1\leq j_2<i_2\leq\ldots\leq
j_{k-1}<i_{k-1}\leq j_k<i_k \quad with\quad  |J|=|I|+h-n.
\]
Now let $\mu=(j_k-k\geq \ldots\geq j_1-1\geq 0)$.
Then:
\[
i_k-k>j_k-k\geq i_{k-1}-(k-1)>j_{k-1}-(k-1)\geq\ldots\geq i_1-1>j_1-1\geq 0,
\]
i.e.
\[
\lambda_1-1\geq\mu_1\geq \lambda_2-1\geq\ldots\geq \lambda_k-1\geq \mu_k\geq 0
\]
Moreover:
\[
|\umu|=|J|-{k(k-1)\over 2}=|I|-{k(k-1)\over 2}+h-n=|\ulamb|+h-n,
\]
proving the claim.

\qed

It follows that the algebra of quantum Schubert calculus is described by the
intersection pair $(\wkM^\vee_n[q],D_q)$. In particular:
\begin{corol}
Any $\sigma_\ulamb$ can be expressed, via the quantum cup product, as a
polynomial in the $\sigma_i$'s with no $q$-correction.
\end{corol}

\proof
\begin{eqnarray*}
&{ }&\sigma_\ulamb\cap
[G]=\Phi_q^([\ikformep])=\Phi_q(\Delta_\ulamb(D)[\kformep])=
\Delta_\ulamb(\sigma)\cap [G].
\end{eqnarray*}

\subsection{Examples}

Here are some computational examples to better show how our algebraic machinery
works.
\subsubsection{Example}
This is the example of the introduction.
\begin{eqnarray*}
&{ }&(\sigma_1\cup_q\sigma_{2,1})\cap
[G]=\sigma_1\cap_q(\sigma_{2,1}\cap[G])=
\sigma_1\cap_q[\epsilon^2\wedge\epsilon^3]=
\Phi_q^{-1}(D_1[\epsilon^2\wedge\epsilon^3])=\\&{ }&
=\Phi_q^{-1}([\epsilon^3\wedge\epsilon^4] -
q[\epsilon^1\wedge\epsilon^2])=\sigma_{2,2}+q\sigma_0.
\end{eqnarray*}

\subsubsection{Example} The iterate of $D_1:\wM\lra\wM$ satisfies the following
formula, which can be checked by induction:
\[
D_1^m(\alpha\wedge\beta)=\sum_{i=0}^m{m\choose i}D_i\alpha\wedge D_{m-i}\beta.
\]
Therefore, on $\bigwedge^2M$:
\begin{eqnarray*}
D_1^4(\epsilon^1\wedge\epsilon^2)+
D_4M&=&\epsilon^5\wedge\epsilon^2+4\epsilon^4\wedge\epsilon^3+
6\epsilon^3\wedge\epsilon^4+4\epsilon^2\wedge\epsilon^5+
\epsilon^1\wedge\epsilon^6=\\
&=&2\epsilon^3\wedge\epsilon^4+3\epsilon^2\wedge\epsilon^5
+\epsilon^1\wedge\epsilon^6.
\end{eqnarray*}
Reading the result in $\bigwedge^2M[D_4]$ one has:
\[
2\epsilon^3\wedge\epsilon^4+3\epsilon^2\wedge\epsilon^5
+\epsilon^1\wedge\epsilon^6=2\epsilon^3\wedge\epsilon^4+2q(\epsilon^2\wedge\epsilon^1)
\]
having set $q:=D_4$. The corresponding expression in quantum Schubert calculus
of $G_1({\PP_{\bf C}}^3)$ is:
\[
\sigma_1\cup_q\sigma_1\cup_q\sigma_1\cup_q\sigma_1=2\sigma_{2,2}+2q\sigma_{0},
\]
(the first member is meaningful by associativity, a non trivial result!)
Hence,
there are
$2$ lines meeting
$4$ lines $\ell_1,\ldots,\ell_4$ in general position in $\PP^3$ and 
$2$ rational
maps $f:(\PP^1,P_1,P_2,P_3,P_4)\lra G_1(\PP^3)$ ($P_i$'s fixed and pairwise
distincts) of degree
$1$ such that $f(P_i)\cap\ell_i\neq \emptyset$.


\end{document}